\documentstyle{amlts}
\begin{document}
\annalsline{157}{2003}
\received{May 22, 2001}
\startingpage{257}
\def\bye{\end{document}}
 \font\tenrm=cmr10

%--------------- Author macros ---------------
%for Bbb in amstex
\catcode`\@=11
\font\twelvemsb=msbm10 scaled 1100
\font\tenmsb=msbm10
%\font\ninemsb=msbm7 scaled 1100%msbm9
\font\ninemsb=msbm10 scaled 800
\newfam\msbfam
\textfont\msbfam=\twelvemsb  \scriptfont\msbfam=\ninemsb
  \scriptscriptfont\msbfam=\ninemsb
\def\msb@{\hexnumber@\msbfam}
\def\Bbb{\relax\ifmmode\let\next\Bbb@\else
 \def\next{\errmessage{Use \string\Bbb\space only in math
mode}}\fi\next}
\def\Bbb@#1{{\Bbb@@{#1}}}
\def\Bbb@@#1{\fam\msbfam#1}
\catcode`\@=12

 \catcode`\@=11
\font\twelveeuf=eufm10 scaled 1100
\font\teneuf=eufm10
\font\nineeuf=eufm7 scaled 1100%eufm9
\newfam\euffam
\textfont\euffam=\twelveeuf  \scriptfont\euffam=\teneuf
  \scriptscriptfont\euffam=\nineeuf
\def\euf@{\hexnumber@\euffam}
\def\frak{\relax\ifmmode\let\next\frak@\else
 \def\next{\errmessage{Use \string\frak\space only in math
mode}}\fi\next}
\def\frak@#1{{\frak@@{#1}}}
\def\frak@@#1{\fam\euffam#1}
\catcode`\@=12
%-------------- Author entries --------------------

%\intro %(Optional, Introduction)
 \newcommand{\sprod}[2]{#1(#2)}
\newcommand{\subspace}{\mathrel{{\subset}{\subset}}}

\newcommand{\R}{{\Bbb R}}
\newcommand{\N}{{\Bbb N}}
\newcommand{\dist}{{\rm dist}}
\newcommand{\diam}{{\rm diam}}
\newcommand{\lspan}{{\rm span}}
\newcommand{\clspan}{\overline{\lspan}}
\newcommand{\Tan}{{\rm Tan}}
\newcommand{\Lip}{{\rm Lip}}
\newcommand{\Lin}{{\rm Lin}}
\newcommand{\ls}{{\cal L}}

\newcommand{\aint}{{-}\hskip-.9em\int}
\newcommand{\Aint}{{-}\hskip-1.05em\int}

\title{On Fr\'echet differentiability
of\\ Lipschitz maps between Banach spaces}
\shorttitle{Fr\'echet differentiability
of Lipschitz maps}  

\twoauthors{Joram Lindenstrauss}{David Preiss}
%\authors{}% Separate each author with a comma and a space.
\institutions{Einstein Institute of Mathematics, The Hebrew University of Jerusalem,
Jerusalem, Israel\\
{\eightpoint {\it E-mail address\/}:  joram@math.huji.ac.il}
\\ \vglue6pt
University College London, London, United Kingdom\\
{\eightpoint {\it E-mail address\/}: d.preiss@ucl.ac.uk
}}
%-------------- Article Text--------------------

\centerline{\bf Abstract}
\vglue12pt

A well-known open question is whether every
countable collection of\break Lipschitz functions
on a Banach space $X$ with separable dual has a common
point of Fr\'echet differentiability.
We show that the answer is positive
for some infinite-dimensional $X$.
Previously, even for collections consisting of two functions
this has been known for finite-dimensional $X$
only (although for one function the answer is
known to be affirmative in full generality).
Our aims are achieved by introducing a new class of
null sets in Banach spaces (called $\Gamma$-null sets),
whose definition involves both the notions of
category and measure, and showing that the required
differentiability holds almost everywhere with respect to it.
We even obtain
existence of Fr\'echet derivatives of Lipschitz
functions between certain infinite-dimensional Banach spaces;
no such results have been known previously.

Our main result states that
a Lipschitz map between separable Banach spaces
is Fr\'echet differentiable
$\Gamma$-almost everywhere provided that it is
regularly G\^ateaux differentiable
$\Gamma$-almost everywhere and the G\^ateaux derivatives
stay within a norm separable space of operators.
It is easy to see that\break Lipschitz maps of $X$
to spaces with the Radon-Nikod\'ym property are
G\^ateaux differentiable $\Gamma$-almost everywhere.
Moreover, G\^ateaux differentiability implies regular
G\^ateaux differentiability with exception of another
kind of negligible sets, so-called $\sigma$-porous sets.
The answer to the question is therefore positive in every space
in which every $\sigma$-porous set is $\Gamma$-null. We show that
this holds for $C(K)$ with $K$ countable compact,
the Tsirelson space and for all subspaces of $c_0$,
but that it fails for Hilbert~spaces.

\section{Introduction}\label{sec:I}

One of the main aims of this paper is to show that infinite-dimensional Banach spaces may have the property that any countable
collection of real-valued Lipschitz functions defined on them has
a common point of Fr\'echet differentiability. Previously, this
has not been known even for collections consisting of two such
functions. Our aims are achieved by introducing a new class of
null sets in Banach spaces and proving results on
differentiability almost everywhere with respect to it. The
definition of these null sets involves both the notions of
category and measure. This new concept even enables the proof of
the existence of Fr\'echet derivatives of Lipschitz
functions between certain infinite-dimensional Banach spaces.
No such results have been known previously.

Before we describe this new class of null sets and the new results, we
present briefly some background material (more details and additional
references can be found in~\cite{BL}).

There are two basic notions of differentiability for functions $f$
defined on an open set in a Banach space $X$ into a Banach space $Y$.
The function $f$ is said to be G\^ateaux differentiable at $x_0$ if
there is a bounded linear operator $T$ from $X$ to $Y$ so that for
every $u\in X$,
$$\lim_{t\to 0} \frac{f(x_0+tu)-f(x_0)}{t} = Tu.$$
The operator $T$ is called the G\^ateaux derivative of $f$ at $x_0$
and is denoted by $D_f(x_0)$.

If for some fixed $u$ the limit
$$f'(x_0,u) = \lim_{t\to 0} \frac{f(x_0+tu)-f(x_0)}{t}$$
exists, we say that $f$ has a directional derivative at $x_0$ in the
direction $u$. Thus $f$ is G\^ateaux differentiable at $x_0$ if and only if
all the directional derivatives $f'(x_0,u)$ exist and they form a
bounded linear operator of $u$. Note that in our notation we have in
this case $f'(x_0,u) = D_f(x_0) u$.

If the limit in the definition of G\^ateaux derivative exists
uniformly in $u$ on the unit sphere of $X$, we say that $f$ is
Fr\'echet differentiable at $x_0$ and $T$ is the Fr\'echet derivative
of $f$ at $x_0$. Equivalently, $f$ is Fr\'echet differentiable at
$x_0$ if there is a bounded linear operator $T$ such that
$$f(x_0+u)=f(x_0)+Tu+o(\|u\|)
\textrm{ as } \|u\| \to 0.$$

It is trivial that if $f$ is Lipschitz and $\dim(X)<\infty$ then the
notion of G\^ateaux differentiability and Fr\'echet differentiability
coincide. The situation is known to be completely different if
$\dim(X)=\infty$. In this case there\break are reasonably satisfactory
results on the existence of G\^ateaux derivatives of Lipschitz
functions, while results on existence of Fr\'echet derivatives are
rare and usually very hard to prove. On the other hand, in many
applications it is important to have Fr\'echet derivatives of $f$,
since they provide genuine local linear approximation to $f$,
unlike the much weaker G\^ateaux derivatives.

Before we proceed we mention that we shall always assume the domain
space to be separable and therefore also $Y$ can be assumed
to be separable.

We state now the main existence theorem for G\^ateaux derivatives.
This is a direct and quite simple generalization of Rademacher's
theorem to infinite-dimensional spaces. But first we recall the
definition of two notions which enter into its statement.

A Banach space $Y$ is said to {\it have the Radon-Nikod{\rm \'{\it y}}m property}
(RNP) if every Lipschitz function $f:\R\to Y$ is differentiable
almost everywhere (or equivalently every such $f$ has a point of
differentiability).

A Borel set $A$ in $X$ is said to be {\it Gauss null\/} if
$\mu(A)=0$ for every
nondegenerate (i.e.\ not supported on a proper closed hyperplane) Gaussian
measure $\mu$ on $X$.
There is also a related notion of Haar null sets which will not be
used in this paper. We just mention, for the sake of orientation, that
the class of Gauss null sets forms a proper subset of the class of
Haar null sets.

\proclaimtitle{\cite{chr}, \cite{man}, \cite{ar}}
\proclaim{Theorem}\label{Gat-orig}
Let $X$ be separable and $Y$ have the {\rm RNP.}
Then every Lipschitz function
from an open set $G$ in $X$ into $Y$
is G{\rm \^{\it a}}teaux differentiable outside a
Gauss null set.
\endproclaim 

In view of the definition of the RNP, the assumption on $Y$ in
Theorem~\ref{Gat-orig} is necessary. Easy and well-known examples
show that Theorem~\ref{Gat-orig} fails badly if we want Fr\'echet
derivatives. For example, the map $f:\ell_2\to\ell_2$ defined by
$f(x_1,x_2,\dots)=(|x_1|,|x_2|,\dots)$ is nowhere Fr\'echet
differentiable.

In the study of Fr\'echet differentiability there is another notion
of smallness of sets
which enters naturally in many contexts. A set $A$ in a Banach space
$X$ (and even in a general metric space) is called {\it porous} if
there is a $0<c<1$ so that for every $x\in A$ there are
$\{y_n\}_{n=1}^\infty\subset X$ with $y_n\to x$ and so that
$B(y_n,c\,\dist(y_n,x))\cap A=\emptyset$ for every $n$. (We denote by
$B(z,r)$ the closed ball with center~$z$ and radius~$r$.)
An important reason for the connections between porous sets and Fr\'echet
differentiability is the trivial remark that if $A$ is porous in a
Banach space $X$ then the Lipschitz function $f(x)=\dist(x,A)$ is not
Fr\'echet differentiable at any point of $A$. Indeed, the only
possible value for the (even G\^ateaux) derivative of $f$ at $x\in A$
is zero. But with $y_n$ and $c$ as above,
$f(y_n) \ge c\,\dist(y_n,x)$ is not $o(\dist(y_n,x))$ as $n\to\infty$.
A set $A$ is called {\it $\sigma$-porous} if it can be represented as a
union
$A=\bigcup_{n=1}^\infty A_n$ of countably many porous sets (the
porosity constant $c_n$ may vary with $n$).

If $U$ is a subspace of a
Banach space $X$, then the set $A$ will be called {\it porous in the
direction $U$} if there is a $0<c<1$ so that for every $x\in A$ and
$\varepsilon > 0$ there is a $u\in U$ with $\|u\| < \varepsilon$
and so that $B(x+u,c\|u\|)\cap A=\emptyset$. A set $A$
in a Banach space $X$ is called {\it directionally porous} if there is a
$0<c<1$ so that for every $x\in A$ there is a $u=u(x)$ with $\|u\|=1$
and a sequence $\lambda_n\searrow 0$ so that
$B(x+\lambda_nu,c\lambda_n)\cap A=\emptyset$ for every $n$. The
notions of {\it $\sigma$-porous sets in the direction $U$} or
{\it $\sigma$-directionally porous} sets are defined in an obvious way.

In finite-dimensional spaces a simple compactness argument shows that
the notions of porous and directionally porous sets
coincide. As it will become presently clear, this is not the case if
$\dim(X)=\infty$. In finite-dimensional spaces porous sets are small
from the point of view of measure (they are of Lebesgue measure zero
by Lebesgue's density theorem) as well as category (they are obviously
of the first category). In infinite-dimensional spaces only the first
category statement remains valid.

As is well known, the easiest class of functions to handle in
differentiation theory are convex continuous real-valued functions
$f:X\to \R$. In~\cite{pz1} it is proved that if $X^\star$ is separable
then any convex continuous
$f:X\to\R$ is Fr\'echet differentiable outside a $\sigma$-porous set.
In separable spaces with $X^\star$ nonseparable it is known~\cite{lw}
that there are convex continuous functions (even equivalent norms)
which are nowhere Fr\'echet differentiable. It is shown in~\cite{mm}
and~\cite{mat} that in every infinite-dimensional super-reflexive
space $X$, and in particular in $\ell_2$, there is an equivalent norm
which is Fr\'echet differentiable only on a Gauss null set. It follows
that such spaces $X$ can be decomposed into the union of
two Borel sets $A\cup B$ with $A$ $\sigma$-porous and $B$ Gauss null.
Such a decomposition was proved earlier and directly for every
separable infinite-dimensional Banach space $X$
(see~\cite{pt}). Note that if $A$ is a directionally porous set in a
Banach space then, by an argument used already above, the Lipschitz
function $f(x)=\dist(x,A)$ is not even G\^ateaux differentiable at any
point $x\in A$ and thus by Theorem~\ref{Gat-orig} the set $A$ is Gauss
null. 

The new null sets (called $\Gamma$-null sets) will be introduced in the next
section. There we prove some simple facts concerning these null sets and
in particular that Theorem~\ref {Gat-orig} also holds if we require
the exceptional set (i.e.\ the set of non-G\^ateaux differentiability)
to be $\Gamma$-null.

The main result on Fr\'echet differentiability in the context of
$\Gamma$-null sets is proved in Section~\ref{sec:F}. From this result
it follows in particular that if every $\sigma$-porous set in $X$ is
$\Gamma$-null then any Lipschitz $f:X\to Y$ with $Y$ having the RNP
whose set of G\^ateaux derivatives $\{D_f(x)\}$ is separable is
Fr\'echet differentiable $\Gamma$-almost everywhere. From the main
result it follows also that convex continuous functions on any
space $X$ with $X^\star$ separable are Fr\'echet differentiable
$\Gamma$-almost everywhere. In particular, if $X^\star$ is separable,
$f:X\to\R$ is convex and continuous and $g:X\to Y$ is Lipschitz with
$Y$ having the RNP then there is a point $x$ (actually
$\Gamma$-almost any point) at which $f$ is Fr\'echet differentiable
and $g$ is G\^ateaux differentiable. This information on existence of
such an~$x$ cannot be deduced from the previously known results. It is
also clear from what was said above that the $\Gamma$-null sets and
Gauss null sets form completely different $\sigma$-ideals in general
(the space $X$ can be decomposed into disjoint Borel sets $A_0\cup
B_0$ with $A_0$ Gauss null and $B_0$ $\Gamma$-null, at least when $X$
is infinite-dimensional and super-reflexive).

In Section~\ref{sec:S} we prove that for $X=c_0$ or more generally
$X=C(K)$ with $K$ countable compact and for some closely related
spaces that every $\sigma$-porous set in them is indeed $\Gamma$-null.
Thus combined with the main result of Section~\ref{sec:F} we get a
general result on existence of points of Fr\'echet differentiability
for Lipschitz maps $f:X\to Y$ where $X$ is as above and $Y$ has the RNP.
This is the first result on existence of points of Fr\'echet
differentiability 
for general Lipschitz mappings for certain pairs of infinite-dimensional spaces. Actually, the only previously known general
result on existence of points of Fr\'echet differentiability of
Lipschitz maps with infinite-dimensional domain
dealt with maps whose range is the real line~\cite{p} and~\cite{lp}.

Unfortunately, the class of spaces in which $\sigma$-porous sets are
$\Gamma$-null does not include the Hilbert space $\ell_2$ or more
generally $\ell_p$, $1<p<\infty$. The reason for this is an example
in~\cite{pt} which shows that for these spaces the mean value theorem
for Fr\'echet derivatives fails while a result in Section~\ref{sec:M}
shows that in the sense of $\Gamma$-almost everywhere the mean value
theorem for Fr\'echet derivatives holds. All this is explained in
detail in Section~\ref{sec:M}.

The paper concludes in Section~\ref{sec:R} with some
comments and open problems.

\section{$\Gamma$-null sets}\label{sec:G}

Let $T=[0,1]^{\N}$ be endowed
 with the product topology and product Lebesgue measure
$\mu$.
Let $\Gamma(X)$ be the space of continuous mappings\break
$\gamma: T\to X$ having continuous partial derivatives
$D_j\gamma$ (with one-sided derivatives at points
where the $j$-th coordinate is~$0$ or~$1$). The elements
of $\Gamma(X)$ will be called surfaces.
For finitely supported $s\in\ell_\infty$ we
also use the notation
$\gamma'(t)(s)=\sum_{j=1}^\infty s_j D_j\gamma(t)$.
We equip $\Gamma(X)$
by the topology generated by the semi-norms
$\|\gamma\|_0^{\phantom{|}}=\sup_{t\in T}\|\gamma(t)\|$
and
$\|\gamma\|_k=\sup_{t\in T}\|D_k\gamma(t)\|$.
Equivalently, this topology may be defined by using
the semi-norms
$\|\gamma\|_{\le k}=\max_{0\le j\le k}\|\gamma\|_j$.
The space $\Gamma(X)$ with this topology is a Fr\'echet space;
in particular, it is Polish (metrizable by a complete separable
metric).

We will often use the simple observation that for every
$\gamma\in\Gamma(X)$, $m\in\N$ and $\varepsilon>0$ there
is $n\in\N$ so that for every $t\in T$
the surface
$$\gamma^{n,t}(s)
=\gamma(s_1,\dots,s_n,t_{n+1},t_{n+2},\dots)$$
satisfies
$$\|\gamma^{n,t}-\gamma\|_{\le m}<\varepsilon.$$
This follows immediately from the uniform continuity of
$\gamma$ and its partial derivatives.
We let $\Gamma_k(X)$ be the space of those $\gamma\in\Gamma(X)$
that depend on the first $k$
coordinates of $T$
and note that by the observation above
$\bigcup_{k=1}^\infty\Gamma_k(X)$ is
dense in $\Gamma(X)$.

The \emph{tangent space} $\Tan(\gamma,t)$  of $\gamma$ at
a point $t\in T$ is defined to be
the closed linear span in $X$ of the vectors
$\{D_k\gamma(t)\}_{k=1}^\infty$.

\numbereddemo{Definition}  \label{def:null}\rm
A Borel set $N\subset X$ will be called $\Gamma$-null if\break
$\mu\{t\in T: \gamma(t)\in N\}=0$ for residually many
$\gamma\in\Gamma(X)$. A possibly non-Borel set $A\subset X$
will be called $\Gamma$-null if it is contained in a Borel
$\Gamma$-null set.
\enddemo

Sometimes, we will also consider $T$ as a subset of $\ell_\infty$.
For example, for $s,t\in T$ we use the notation
$\|s-t\|=\sup_{j\in\N} |s_j-t_j|$.
We also use the notation
$Q_k(t,r)=\{s\in T:\max_{1\le j\le k} |s_j-t_j|\le r\}$.

\proclaim{Lemma}\label{tan-lem}
Let $\{u_j\}_{j=1}^n\subset X$ and $\varepsilon>0$. Then the set of those
$\gamma\in\Gamma(X)$ for which there are $k\in\N$ and
$c>0$ such that
$$\max_{1\le j \le n}
\sup_{t\in T}\|c D_{k+j}\gamma(t) - u_j\|<\varepsilon $$
is dense and open in $\Gamma(X)$.
\endproclaim

\demo{Proof} By the definition of the topology of $\Gamma(X)$ it is clear that
this set is open. 
To see that it is dense it suffices to show that its closure
contains $\Gamma_k(X)$ for every $k$.
Let $\gamma_0\in\Gamma_k(X)$,
$\eta >0$ and consider the surface
$\gamma(t)=\gamma_0(t)+ \eta \sum_{j=1}^n t_{k+j} u_j$.
Then 
$\|\gamma-\gamma_0\|_0 \le
n \eta \max_{1\le j \le n} \|u_j\|$,
$\|\gamma-\gamma_0\|_l=0$ if $l\le k$ or $l>k+n$,
$\|\gamma-\gamma_0\|_{k+j}=\eta\|u_j\|$ and
$c D_{k+j}\gamma(t)= u_j$ for $1\le j \le k$
and $c = 1/\eta$.
\enddemo

\proclaim{{C}orollary}\label{tan-res}
If $X$ is separable{\rm ,} then  residually many
$\gamma\in\Gamma(X)$ have the property that
$\Tan(\gamma,t)=X$ for every $t\in T$.
\endproclaim 

\demo{Proof} By Lemma~\ref{tan-lem} (with $n=1$) we get that for every $u\in X$
the set of those $\gamma\in\Gamma(X)$ such that
$\dist(u,\Tan(\gamma,t))<\varepsilon$ for every $t\in T$ is open and
dense in $\Gamma(X)$. The desired result follows now from the
separability of $X$.
\enddemo

We show next that the class of $\Gamma$-null sets in a finite-dimensional space $X$ coincides with the class of sets of
Lebesgue measure zero (just as for Gauss and Haar null sets).

\proclaim{Theorem}\label{fd}
In finite\/{\rm -}\/dimensional spaces{\rm ,} $\Gamma$\/{\rm -}\/null sets coincide with
Le\-besgue null sets.
\endproclaim 

\demo{Proof} Let $u_1,\dots,u_n\in X$ be a basis for $X$, let $E\subset X$ be
a Borel set and denote by $|E|$ its Lebesgue measure.

If $|E| > 0$, define $\gamma_0: T\to X$ by
$\gamma_0(t)=u+\sum_{j=1}^n t_j u_j$, where $u\in X$ is
chosen so that $|E\cap\gamma_0(T)|>0$.
If $\|\gamma-\gamma_0\|_{\le n}$ is sufficiently small,
then for every $s=(s_1,s_2,\dots)\in T$
the mappings $\gamma_s(t_1,\dots,t_n)
=\gamma(t_1,\dots,t_n,s_1,s_2,\dots)$
are diffeomorphisms of $[0,1]^n$ onto subsets of $X$ which meet $E$
in a set of measure at least $|E\cap\gamma_0(T)|/2$. Hence,
for every $s$,
$|\gamma_s^{-1}(E)| \ge c_1|E\cap\gamma_0(T)|$,
for a suitable positive constant
$c_1$. Hence
$$\mu(\gamma^{-1}(E))=\int_T |\gamma_s^{-1}(E)|\,d\mu(s)
\ge c_1|E\cap\gamma_0(T)|$$
and we infer that $E$ is not $\Gamma$-null.

If $|E|=0$, we use
Lemma~\ref{tan-lem} with a sufficiently small $\varepsilon>0$
to find a dense open set of such surfaces
$\gamma\in\Gamma(X)$ for which there are $k\in\N$ and
$c>0$ such that
$$\max_{1\le j \le n}
\sup_{t\in T}\| D_{k+j} c\gamma(t) - u_j\|<\varepsilon.$$
Then the mappings 
$c\gamma_s(t_1,\dots,t_n)=
c\gamma(s_1,\dots,s_k,t_1,\dots,t_n,s_{k+1},s_{k+2},\dots)$
are, for every $s\in T$,
diffeomorphisms of $[0,1]^n$ onto a subset of $X$.
The same is therefore true for the mappings
$\gamma_s$, $s\in T$.
Hence
$|\gamma_s^{-1}(E)|=0$ for every $s$ and hence
$$\mu(\gamma^{-1}(E))=\int_T |\gamma_s^{-1}(E)|\,d\mu(s)=0;$$
i.e.\ $E$ is $\Gamma$-null.
\enddemo

We show next that Theorem~\ref{Gat-orig} remains valid if we replace
in its statement Gauss null sets by $\Gamma$-null sets.

\proclaim{Theorem}\label{gat-ae}
Let $X$ be separable and $Y$ have the {\rm RNP.}
Then every Lipschitz function
from an open set $G$ in $X$ into $Y$
is G{\rm \^{\it a}}teaux differentiable outside a
$\Gamma$\/{\rm -}\/null set.
\endproclaim 

\demo{Proof} We remark first that the set of points at which $f$ fails to be
G\^ateaux differentiable is a Borel set.
We recall next that
Rademacher's theorem holds also for Lipschitz maps
from $\R^k$ to a space $Y$ having the RNP
(see e.g.\ \cite[Prop.~6.41]{BL}).
Consider now an arbitrary surface $\gamma$. By using
Fubini's theorem, we get that for almost every $t\in T$
for which $\gamma(t)\in G$ the
mapping 
$$(s_1,\dots,s_k)\to f(\gamma(s_1,\dots,s_k,t_{k+1},\dots))$$
is
differentiable at $s=(t_1,\dots,t_k)$. Since $f$ is Lipschitz,
it follows
that, for almost all $t$, $f$ has directional derivatives
for all vectors 
$$v \in \clspan\{D_j(\gamma)(t)\}_{j=1}^\infty=\Tan(\gamma,t)$$
at $u=\gamma(t)$ and that these directional derivatives depend
linearly on $v$ (see e.g.\ \cite[Lemma 6.40]{BL}).
In particular, for every surface $\gamma$ from the residual set
obtained in
Corollary~\ref{tan-res}
$f$ is G\^ateaux differentiable at $u=\gamma(t)$
for almost all $t\in T$ for which $\gamma(t)\in G$.
This proves the theorem.
\enddemo

\section{Fr\'echet differentiability}\label{sec:F}

In this section we prove the main criterion for Fr\'echet differentiability
of Lipschitz functions in terms of $\Gamma$-null sets. But first we
have to introduce the following simple notion.

\numbereddemo{Definition}  \label{def:reg}
Suppose that $f$ is a map from (an open set in) $X$ to $Y$.
We say that a point $x$
is a \emph{regular point of $f$} if for every $v\in X$ for which
$f'(x,v)$ exists,
$$\lim_{t\to 0} \frac{f(x+tu+tv)-f(x+tu)}{t}=f'(x,v)$$
uniformly for $\|u\|\le 1$.
\enddemo

Note that in the definition above
it is enough to take the limit for $t\searrow 0$
only, since we may replace $v$ by $-v$.

\proclaim{Proposition} \label{conv-reg}
For a convex continuous function $f:X\to\R$ every
point $x$ is a regular point of $f$.
\endproclaim 

\demo{Proof} Given $x\in X$, $v\in X$ and $\varepsilon>0$, find $r>0$ such that
$$|(f(x+tv)-f(x))/t -f'(x,v)|<\varepsilon$$ for $0<|t|<r$
and such that $f$ is Lipschitz on $B(x,2r(1+\|v\|))$ with constant
$K$. If $\|u\|\le 1$ and $0<t<\min(r,\varepsilon r/2K)$,
then
\begin{eqnarray*}
(f(x+tu+tv)-f(x+tu))/t &\le& (f(x+tu+rv)-f(x+tu))/r\\
&\le&
(f(x+rv)-f(x))/r+ 2Kt\|u\|/r \\
&<& f'(x,v)+2\varepsilon\\
\noalign{\noindent and}
(f(x+tu+tv)-f(x+tu))/t &\ge& (f(x+tu)-f(x+tu-rv))/r\\
&\ge&
(f(x)-f(x-rv))/r- 2Kt\|u\|/r \\
&>& f'(x,v)-2\varepsilon.\\
\noalign{\vskip-36pt}
\end{eqnarray*}
\enddemo

{\it Remark}.
It is well known and as easy to prove as the  statement above
that convex functions satisfy a stronger condition of regularity
(sometimes called Clarke regularity),
namely that
$$\lim_{z\to x,\, t\to 0} \frac{f(z+tv)-f(z)}{t}=f'(x,v)$$
whenever $f'(x,v)$ exists.
We do not use here this stronger regularity concept since while
every point of Fr\'echet differentiability of $f$ is a point of
regularity of $f$ in our sense, this no longer holds for the stronger
regularity notion; this is immediate by considering an indefinite
integral of the characteristic function of a set $E\subset\R$
such that both $E$ and its complement have positive measure
in every interval.
Therefore the stronger form of regularity cannot be used in
proving existence of points of differentiability for Lipschitz maps
(which is our purpose here).

\proclaim{Proposition} \label{por-prp}
Let $f$ be a Lipschitz map from an open subset $G$ of a separable
Banach space $X$ to a separable Banach space $Y$.
Then the set of irregular points of $f$ is $\sigma$\/{\rm -}\/porous.
\endproclaim 

\demo{Proof} For $p,q\in\N$, $v$ from a countable dense subset of $X$
and $w$ from a countable dense subset of $Y$
let
$E_{p,q,v,w}$ be the set of those $x\in X$ such that
$\|f(x+tv)-f(x)-tw\|\le |t|/p$ for $|t|<1/q$, and
$$\limsup_{t\to 0}\sup_{\|u\|\le 1}
\|\frac{f(x+tu+tv)-f(x+tu)}{t}-w\|>2/p.$$

Whenever $x\in E_{p,q,v,w}$, there are
arbitrarily small $|t|<1/q$
such that for some $u$ with $\|u\|\le 1$,
$$\|f(x+tu+tv)-f(x+tu)-tw\|>2|t|/p.$$
If $\|y-(x+tu)\|<|t|/2p\Lip(f)$,
then
$$\|f(y+tv)-f(y)-tw\|\ge \|f(x+tu+tv)-f(x+tu)-tw\|
-|t|/p>|t|/p$$ 
and hence $y\notin E_{p,q,v,w}$. This proves that
$E_{p,q,v,w}$ is $1/2p\Lip(f)$ porous.
Since every irregular point of $f$ belongs to some
$E_{p,q,v,w}$ the result follows.
\enddemo

The next lemma is a direct consequence of the definition of regularity.
It will make the use of the regularity assumption more convenient
in subsequent arguments.

\proclaim{Lemma}\label{reg-lem}
Suppose that $f$ is Lipschitz on a neighborhood of $x$
and that{\rm ,} at $x${\rm ,} it
is regular and differentiable in the direction of a
finite\/{\rm -}\/dimensional subspace
$V$ of $X$. Then for every $C$ and $\varepsilon>0$ there is
a $\delta>0$
such that
$$\|f(x+v+u)-f(x+v)\|\ge\|f(x+u)-f(x)\|-\varepsilon\|u\|$$
whenever $\|u\|\le\delta$, $v\in V$ and $\|v\|\le C\|u\|$.
\endproclaim

\demo{Proof} Let $r>0$ be such that $f$ is Lipschitz on  $B(x,r)$.
Let $S$ be a finite subset of $\{v\in V:\|v\|\le C\}$
such that for every $w$ in this set
there is $v\in S$
such that $\|w-v\|<\varepsilon/6\Lip(f)$.
By the definition of regularity, there is $0<\delta < r/(1+C)$
such that
\begin{equation}\label{preg.1}
\|f(x+t\hat u+t\hat v)-f(x+t\hat u)-tf'(x,\hat v)\|\le\varepsilon t/3
\end{equation}
whenever $0\le t\le\delta$, $\|\hat u\|\le 1$ and $\hat v\in S$.

Suppose that $\|u\|\le\delta$, $v\in V$ and $\|v\|\le C\|u\|$.
By using (\ref{preg.1}) with $t=\|u\|$,
$\hat u=u/\|u\|$ and $\hat v\in S$
with $\|\hat v - v/\|u\|\|<\varepsilon/6\Lip(f)$,
we get
\begin{equation}\label{preg.2}
\|f(x+u+t\hat v)-f(x+u)-tf'(x,\hat v)\|\le\varepsilon t/3.
\end{equation}
Similarly by using (\ref{preg.1}) with
the same $t$ and $\hat v$ but with $\hat u=0$, we get
\begin{equation}\label{preg.3}
\|f(x+t\hat v)-f(x)-tf'(x,\hat v)\|\le\varepsilon t/3.
\end{equation}
Hence by the triangle inequality, (\ref{preg.2}) and (\ref{preg.3})
we deduce that
$$\|(f(x+t\hat v+u)-f(x+t\hat v))-(f(x+u)-f(x))\|
\le 2\varepsilon t/3.$$
Recalling that $t=\|u\|$ we thus get
\begin{eqnarray*}
\lefteqn{\|f(x+v+u)-f(x+v)\|}\phantom{xxxx}\\
&\ge&
\|(f(x+t\hat v+u)-f(x+t\hat v))\|-2\Lip(f)\|v-t\hat v\|\\
&\ge&
\|f(x+u)-f(x)\|-2\varepsilon t/3 -2\Lip(f)t(\varepsilon/6\Lip(f))\\
&=&
\|f(x+u)-f(x)\|-\varepsilon\|u\|.\\
\noalign{\vskip-36pt}
\end{eqnarray*}
\enddemo
\vglue12pt
Our next lemma examines the consequence
of having a Lipschitz function which is not
Fr\'echet differentiable at a regular point:
It shows that if a surface contains such a point $x$
then, after a suitable small perturbation of the surface,
the derivative near $x$, in the mean, is not close to its value
at $x$. Moreover, this property persists for all surfaces close
enough to the perturbed surface.
As in the proof
of Lemma~\ref{tan-lem} (and other proofs in \S\ref{sec:G})
we will make here an essential use of the fact that in the context
of $\Gamma$-null sets we work with infinite-dimensional surfaces
$\gamma\in\Gamma(X)$. These surfaces can be well approximated by
$k$-dimensional surfaces in $\Gamma_k(X)$. The surfaces in
$\Gamma_k(X)$ can in turn be approximated by surfaces in
$\Gamma_{k+1}(X)$ and we are quite free to do appropriate constructions
on the $(k+1)$-th coordinate in order to get an approximation with
desired properties.

\proclaim{Lemma}\label{nonfpt-lem}
Let $f:G\to Y$ be a Lipschitz function with $G$ an open set in a
separable Banach space $X$. Let $E$ be a Borel subset of $G$
consisting of points where $f$ is G{\rm \^{\it a}}teaux differentiable
and regular.
Let $\eta>0${\rm ,}
$\tilde\gamma\in\Gamma_k(X)${\rm ,} $t\in T$
so that $x=\tilde\gamma(t)\in E$.

Then there are $0<r<\eta${\rm ,} $\delta>0$ and
$\hat\gamma\in\Gamma_{k+1}(X)$ such
that
$\|\hat\gamma-\tilde\gamma\|_{\le k}<\eta${\rm ,}
$\hat\gamma(s)=\tilde\gamma(s)$ for $s\in T\setminus Q_k(t,r)$
and so that the following holds\/{\rm :} Whenever
$\gamma\in\Gamma(X)$ has the property that
$$\|\gamma(s)-\hat\gamma(s)\|
+\|D_{k+1}(\gamma(s)-\hat\gamma(s))\| <\delta,
\quad s\in Q_k(t,r)$$
then either
$$\mu(Q_k(t,r)\setminus\gamma^{-1}(E))
\ge\alpha\mu(Q_k(t,r))/8\Lip(f)$$
or
$$\int_{Q_k(t,r)\cap\gamma^{-1}(E)}
\|D_f(\gamma(s)) -D_f(x)\|\,
  d\mu(s)\ge\alpha\mu(Q_k(t,r))/2,$$
where  $\alpha=\limsup_{u\to 0}\|f(x+u)-f(x)-D_f(x)u\|/\|u\|$.
\endproclaim

\demo{Proof} We shall assume that $\alpha > \eta>0$ and
will define $r$ and $\delta$ later in the proof.
Let
$0<\zeta<\min(1,\alpha/5)$ be such that
$$(\alpha-5\zeta)\mu(Q_k(t,(1-\zeta)r))\ge
3\alpha\mu(Q_k(t,r))/4$$
whenever $t\in T$
and $r<1$. This can evidently be accomplished even though
for some $t$, $Q_k(t,r)$ is not an entire cube.

Choose a continuously differentiable function
$\omega:\ell_\infty\to [0,1]$
depending on the first $k$ coordinates only
such that $\omega(s)=1$ if
$\|\pi_k s\|\le 1-\zeta$ and
$\omega(s)=0$ if $\|\pi_k s\|\ge 1$
where $\pi_k$ denotes the projection of $\ell_\infty$
on its first $k$ coordinates.
Let $\max(4,\alpha)\le K<\infty$ be such
that
$\|\tilde\gamma(s_1)-\tilde\gamma(s_2)\|\le K\|s_1-s_2\|$,
$\|\omega(s_1)-\omega(s_2)\|\le K\|s_1-s_2\|$,
and the function
$$g(z)=f(z)-f(x)-D_f(x)(z-x)$$
is Lipschitz with constant $K$ on
$B(x,2(1+K^2/\eta)\delta_1)\subset G$ for some
$0<\delta_1<\eta$.

We let $C=4K^2/\eta$ and use Lemma~\ref{reg-lem}
to find a $0<\delta_2<\delta_1<\eta$ such that
\begin{equation}\label{pnf.1}
\|g(x+v+u)-g(x+v)\|
\ge\|g(x+u)\|-\zeta\|u\|
\end{equation}
whenever
$\|u\|\le\delta_2$, $v\in\Tan(\tilde\gamma,t)$
(which is of dimension at most $k$ since
$\tilde\gamma\in\Gamma_k(X)$)
and $\|v\|\le C\|u\|$.
Also, let $\delta_3>0$ be such that
\begin{equation}\label{tgm.1}
\|\tilde\gamma(s)-x-\tilde\gamma'(t)(s-t)\|
\le\zeta\|s-t\|/C
\end{equation}
whenever $s\in T$ and $\|\pi_k(s-t)\|\le \delta_3$.
By the definition of $\alpha$
we can find a $u\in X$ such that
$$0<\|u\|<\min(\delta_2,2K\delta_3/C,\eta^2/2K)$$ and
$\|g(x+u)\|\ge(\alpha-\zeta)\|u\|$.

Whenever $s\in T$ and
$\|\pi_k(s-t)\|< C\|u\|/2K$, we use~(\ref{tgm.1})
(which is applicable since  $C\|u\|/2K<\delta_3$)
and~(\ref{pnf.1}) with
$v=\tilde\gamma'(t)(s-t)$
(this is justified by $\|\tilde\gamma'(t)(s-t)\|
\le K\|\pi_k(s-t)\|<C\|u\|$) to infer that
\begin{eqnarray}
\lefteqn{\|g(\tilde\gamma(s)+u)-g(\tilde\gamma(s))\|}\phantom{xxxx}\label{al
p3}\\
&\ge&
\|g(x+v+u)-g(x+v)\|-2K\|\tilde\gamma(s)-(x+v)\|
\nonumber\\
&\ge&
\|g(x+v+u)-g(x+v)\|-2K\zeta\|s-t\|/C\nonumber\\
&\ge&
\|g(x+u)\|-\zeta\|u\| - \zeta\|u\|\nonumber\\
&\ge&
(\alpha-3\zeta)\|u\|.\nonumber
\end{eqnarray}

We now define $r$ to be $2K\|u\|/\eta$ and put
$\hat\gamma(s) =\tilde\gamma(s)+ s_{k+1}\omega((s-t)/r)u$.
By the choice of $\|u\|$ we have $0<r<\eta$.
Also,
$\|\hat\gamma-\tilde\gamma\|\le \|u\|<\eta$
and
$\|D_j(\hat\gamma-\tilde\gamma)\|\le K\|u\|/r<\eta$
for $1\le j\le k$ and thus
$\|\hat\gamma-\tilde\gamma\|_{\le k}<\eta$
as required.

We show that the statement of the lemma holds with
$\delta=\zeta\|u\|/2K$. Suppose that
$$\|\gamma(s)-\hat\gamma(s)\|
+\|D_{k+1}(\gamma(s)-\hat\gamma(s))\|
<\delta \textrm{ for } s\in Q_k(t,r)$$
and that 
$$\mu(Q_k(t,r)\setminus\gamma^{-1}(E))
< \alpha\mu(Q_k(t,r))/8\Lip(f).$$
Consider any $s\in Q_k(t,(1-\zeta)r)$ and put
$$\overline{s}=s+(1-s_{k+1})e_{k+1} \textrm{ and }
\underline{s}=s-s_{k+1}e_{k+1}$$
where $e_{k+1}$ denotes the $(k+1)$-th unit vector in $\ell_\infty$,
$I_s=[-s_{k+1},1-s_{k+1}]$ and
$J_s=\{\sigma\in I_s: \gamma(s+\sigma e_{k+1})\in E\}$.
Let $u^\star\in Y^\star$ be such that $\|u^\star\|=1$ and
$\sprod{u^\star}{g(\gamma(\overline{s}) )
-g(\gamma(\underline{s}))}
= \|g(\gamma(\overline{s}) ) -g(\gamma(\underline{s}))\|$.
Then
$\|\pi_k(s-t)\|<r=2K\|u\|/\eta = C\|u\|/2K$ and by
using (6) and the inequalities
$\|\gamma(\overline{s}) - \hat\gamma(\overline{s})\|,
\|\gamma(\underline{s}) - \hat\gamma(\underline{s})\|<\delta$
and $\Lip(g)\le K$,
we get
\begin{eqnarray*}
\int_{I_s}\textstyle
\frac{\partial\sprod{u^\star}{g(\gamma(s+\sigma e_{k+1}))}}
 {\partial\sigma}\, d\sigma
&=& \|g(\gamma(\overline{s}) ) -g(\gamma(\underline{s}))\|\\
&\ge&
\|g(\hat\gamma(\overline{s}) ) -g(\hat\gamma(\underline{s}))\|
-\zeta\|u\|\\
&=& \|g(\tilde\gamma(s)+u)-g(\tilde\gamma(s))\| -\zeta\|u\|\\
&\ge&
(\alpha-4\zeta)\|u\|.
\end{eqnarray*}
Note that since $s\in Q_k(t,(1-\zeta)r)$ we have
$\frac{\partial\hat\gamma(s+\sigma e_{k+1})}
  {\partial\sigma} = u$ and
$$\|D_{k+1}(\gamma(s+\sigma e_{k+1})-\hat\gamma(s+\sigma e_{k+1}))\|
<\delta$$ and hence $\|\frac{\partial\gamma(s+\sigma e_{k+1})}
  {\partial\sigma} - u\|<\delta$  for all $\sigma\in I_s$.
Using also that 
$\Lip(g)\le 2\Lip(f)$ and that $K\ge\max(4,\alpha)$, we infer that
\begin{eqnarray*}
\lefteqn{\int_{J_s}\|D_g(\gamma(s+\sigma
e_{k+1}))\|\,d\sigma}\phantom{xxx}\\
&\ge&
\int_{I_s} |\textstyle
\frac{\partial\sprod{u^\star}{g(\gamma(s+\sigma e_{k+1}))}}
 {\partial\sigma}|\, d\sigma/(\|u\|+\delta)
- |I_s\setminus J_s|\Lip(g)\\
&\ge& (\alpha-4\zeta)/(1+\zeta/2K) -
 2|I_s\setminus J_s|\Lip(f)\\
&\ge& (\alpha-4\zeta) (1-\zeta/2K) -
 2|I_s\setminus J_s|\Lip(f)\\
&\ge& \alpha-5\zeta -
 2|I_s\setminus J_s|\Lip(f).
\end{eqnarray*}
Integrating over $s$ we obtain the desired \pagebreak result:
\begin{eqnarray*}
\lefteqn{\int_{Q_k(t,r)\cap \gamma^{-1}(E)}
\|D_g(\gamma(s))\|\,
  d\mu(s)}\phantom{xxx}\\
&\ge&
\int_{Q_k(t,(1-\zeta)r)\cap \gamma^{-1}(E)}
\|D_g(\gamma(s))\|\,
  d\mu(s)\\
&\ge& (\alpha-5\zeta)\mu(Q_k(t,(1-\zeta)r))
-  2\Lip(f)\mu(Q_k(t,r)\setminus \gamma^{-1}(E))\\
&\ge& 3\alpha\mu(Q_k(t,r))/4 - 2\Lip(f)\alpha\mu(Q_k(t,r))/8\Lip(f)\\
&\ge&
\alpha\mu(Q_k(t,r))/2.\\
\noalign{\vskip-36pt}
\end{eqnarray*}
\enddemo

We now extend the notation introduced in the beginning of Section~2
for $\gamma\in\Gamma(X)$ to the more general setting of
$L_1(T,X)$. For $k\in\N$ and $s\in T$, $g\in L_1(T,X)$ we put
$$g^{k,s}(t)=g(t_1,\dots,t_k,s_{k+1},\dots).$$
With this notation, the Fubini's formula says that
for every $k$,
$$\int_T g(t)\,d\mu(t)=\int_T\int_T g^{k,s}(t)\,d\mu(s)\,d\mu(t)
=\int_T\int_T g^{k,s}(t)\,d\mu(t)\,d\mu(s).$$

As in the case of continuous functions, the functions $g^{k,s}$
approximate $g$ for large enough $k$. The precise formulation of this
statement is given in

\proclaim{Lemma}\label{L1.2}
Suppose that $g\in L_1(T,X)$. Then for any $\eta >0$ there is an
$l\in\N$ such that 
$$\mu\{s\in T: \|g^{k,s}-g\|_{L_1} > \eta\}<\eta$$
for $k\ge l$.
\endproclaim

\demo{Proof} Let $\tilde g:T\to X$ be a continuous function depending on the
first $l$ variables only such that
$\|g-\tilde g\|_{L_1}<\eta^2$
and let $k\ge l$. By Fubini's theorem,
$$\int_T \|g^{k,s}-\tilde g\|_{L_1}\,d\mu(s)
=\|g - \tilde g\|_{L_1} < \eta^2$$
so by Chebyshev's inequality
$$\mu\{s\in T:  \|g^{k,s}-\tilde g\|_{L_1} > \eta\} < \eta$$
as desired.
\enddemo

The next lemma is a version of Lebesgue's differentiability theorem for
functions defined on the infinite torus $T$.

\proclaim{Lemma}\label{L1.1}
Suppose that $g\in L_1(T,X)$. Then for every $\kappa>0$ there is an
$l\in\N$ such that for all $k\ge l$ there are $\delta>0$ and
$A\subset T$ with $\mu(A)<\kappa$ such that
$$\int_{Q_k(t,r)}\|g(s)-g(t)\| \,d\mu(s)<\kappa\mu(Q_k(t,r))$$
for every $t\in T\setminus A$ and $0<r<\delta$.
\endproclaim

\demo{Proof} We find it convenient to use the notation
$$\Aint_{Q_k(t,r)} h(s)\,d\mu(s) \quad\textrm{for}\quad
\mu(Q_k(t,r))^{-1}\int_{Q_k(t,r)} h(s)\,d\mu(s).$$

Let $\tilde g:T\to X$ be a continuous function depending on the
first $l$ coordinates so that
$\|g-\tilde g\|_{L_1}<\kappa^2/9$.
If we put
$N=\{t:\|g(t)-\tilde g(t)\|\ge\kappa/3\}$, we get
by Chebyshev's inequality that $\mu(N) < \kappa/3$.

Fix an arbitrary $k\ge l$ and put
$h(s)=\int\|g^{k,\sigma}(s)-\tilde g(s)\|d\mu(\sigma)$
and $S=\{s:h(s)\ge \kappa/3\}$. Then as in the proof
of Lemma~\ref{L1.2} we get from Fubini's theorem that
$\mu(S)<\kappa/3$. 

Let 
$$A_n=\{t\notin S: \Aint_{Q_k(t,r)} h(s) d\mu(s) \ge \kappa/3
\textrm{ for some } 0<r<1/n\}.$$
By Lebesgue's differentiability theorem,
$$\lim_{r\to 0}\Aint_{Q_k(t,r)} h(s) d\mu(s) = h(t)$$
for $\mu$-almost every $t$. Hence $A_n$ is a decreasing sequence
of measurable sets whose intersection has measure zero.
Let $n_0$ be such that $\mu(A_{n_0})< \kappa/3$
and put
$A=N\cup S\cup A_{n_0}$.
Choose $0<\delta<1/n_0$ such that
$\|\tilde g(s)-\tilde g(t)\|<\kappa/3$
whenever $s\in Q_k(t,\delta)$. Then $\mu(A)<\kappa$
and for every $t\in T$ and $0<r<\delta$ we have,
by Fubini's theorem and the fact that $\tilde g$ depends on the first
$k$ coordinates, that
\begin{eqnarray*}
\int_{Q_k(t,r)}\|g(s)-\tilde g(s)\|d\mu(s)
&=&\int_{Q_k(t,r)}\int_T
  \|g^{k,\sigma}(s)-\tilde g(s)\|d\mu(\sigma)d\mu(s)\\
&=&\int_{Q_k(t,r)}h(s)d\mu(s).
\end{eqnarray*}
Hence, for every $t\in T\setminus A$ and $0<r<\delta$,
\begin{eqnarray*}
 \Aint_{Q_k(t,r)}\|g(s)-g(t)\|d\mu(s)
&\le& \Aint_{Q_k(t,r)}\|g(s)-\tilde g(s)\|d\mu(s)
  + \|g(t)-\tilde g(t)\|\\
  && + \Aint_{Q_k(t,r)}\|\tilde g(s)-\tilde g(t)\|d\mu(s)
  \\
  &<& \kappa/3 + \kappa/3 + \kappa/3 = \kappa.\\
\noalign{\vskip-36pt}
\end{eqnarray*}
\enddemo

\vglue12pt
Our next lemma is in the spirit of descriptive set theory
(see e.g.\ \cite{kechris} for background). It is this lemma
which makes it clear why the separability assumption on
$\ls$ in the statement of the main theorem
(Theorem~\ref {dif-gen2} below) is needed.

Before stating the lemma, we list some assumptions and definitions
which enter into its statement. We assume that $X$ and $Y$ are
separable Banach spaces. The function $f$ is a Lipschitz
function from an open set $G\subset X$ into~$Y$. We let
$\ls$ be a norm separable subspace of the space $\Lin(X,Y)$
of bounded linear operators from $X$ to $Y$. We let $E$ be the set
of those points $x$ in $G$ at which $f$ is regular and G\^ateaux
differentiable and $D_f(x)\in\ls$. We denote by $\varphi$ the
characteristic 
function of $E$ (as a subset of $X$) and let
$\psi(x)=D_f(x)$ for $x\in E$ and
$\psi(x)=0$ if $x\notin E$.
We also put $\Phi(\gamma)=\varphi\circ\gamma$ and
$\Psi(\gamma)=\psi\circ\gamma$.

\proclaim{Lemma}\label{meas3-lem}
The set $E$ is a Borel set and the mappings $\varphi:X\to\R${\rm ,}
 $\psi:X\to\ls${\rm ,} $\Phi:\Gamma(X)\to L_1(T)$ and
$\Psi:\Gamma(X)\to L_1(T,\ls)$ are all Borel measurable.
In particular there is a residual subset $H$ of $\Gamma(X)$
such that the restrictions of
$\Phi$ and $\Psi$ to $H$ are continuous.
\endproclaim

\demo{Proof} For $L\in\ls$, $u,v\in X$ and $\sigma,\tau>0$ denote by
$M(L,u,v,\sigma,\tau)$ the set of $x\in G$ such that
$\dist(x,X\setminus G)\ge \tau(\|u\|+\|v\|)$ and
$\|f(x+su+sv)-f(x+sv)-sL(u)\|\le |s|\sigma\|u\|$ whenever $|s|<\tau$.
Clearly each $M(L,u,v,\sigma,\tau)$ is a closed
subset of $X$.

Let ${\cal S}$ be a countable dense subset of $\ls$,
$D$ a dense countable subset of $X$, and
$R$ be the set of positive rational numbers. Then
$$E= \bigcap_{\sigma\in R}\bigcup_{L\in{\cal S}}
    \bigcap_{u\in D}\bigcup_{\tau\in R}
    \bigcap_{v\in D, \|v\|\le 1}
    M(L,u,v,\sigma,\tau)$$
and hence $E$ is Borel and $\varphi$ is Borel measurable.

For every
$L\in\ls$ and $\varrho>0$ we have
$$\{x\in E: \|\psi(x)-L\|\le \varrho\} =E\cap \bigcap_{u\in D}
 \bigcap_{R\ni \sigma>\varrho}\bigcup_{\tau\in R}
 M(L,u,0,\sigma,\tau).$$
Since $E$ is Borel, $\ls$ is separable and $\psi(x)=0$
outside $E$, it follows that also $\psi$ is Borel
measurable.

Since $\psi$ is bounded and Borel measurable, the Borel measurability
of $\Psi$ will be established once we show that for every Borel
measurable bounded $h:X\to\ls$ the mapping
$\Psi_h:\Gamma(X)\to L_1(T,\ls)$ defined by
$\Psi_h(\gamma)=h\circ\gamma$ is Borel measurable.
If $h$ is continuous, then so is $\Psi_h$.  If
$\{h_n\}_{n=1}^\infty$ are uniformly bounded in norm and
$h_n\to h$ pointwise and if all
$\Psi_{h_n}$ are Borel measurable, then
$\Psi_{h_n}$ converge pointwise to $\Psi_h$, so $\Psi_h$ is Borel
measurable. The same argument shows the Borel measurability of $\Phi$.

The last statement in the lemma follows from the general fact that
a Borel measurable mapping between complete separable metrizable
spaces has a continuous restriction to a suitable residual subset.
\enddemo

\demo{{R}emark} 
Without assuming the separability of $\ls$ not only does the proof not
work but the statement is actually false. The Borel image of a complete
separable metric space is again separable. Hence if $\Psi$ is Borel,
the image of $\Gamma(X)$ under $\Psi$ must be separable.
\enddemo

The final lemma before the proof of the main theorem combines much of
the previous lemmas. According to Lemma~\ref{nonfpt-lem},
if $f$ is regular but not Fr\'echet differentiable
at a point of a surface then a suitable small local deformation of the
surface causes a nonsmall perturbation of
the function $t\to D_f(\gamma(t))$. Here we perform a finite number
of such local perturbations (on suitable disjoint neighborhoods
chosen via the Vitali theorem) and get the same effect globally.
The notation is as in the preceding lemmas.

\proclaim{Lemma}\label{nonf-lem}
Suppose that $\varepsilon,\eta>0$ and that
$\tilde\gamma\in\Gamma_l(X)$ is such that the set
$$S=\{t\in\tilde\gamma^{-1}(E):
 \limsup_{u\to 0}
\|f(\tilde\gamma(t)+u)-f(\tilde\gamma(t))
    -D_f(\tilde\gamma(t))u\|/\|u\|>\varepsilon\}$$
has $\mu$ measure greater than $\varepsilon$.
Then there are $k\ge l${\rm ,} $\delta>0$ and
$\hat\gamma\in\Gamma_{k+1}(X)$
such that $\|\hat\gamma-\tilde\gamma\|_{\le k}<\eta$
and
$$\int_T (|\varphi(\gamma(t))-\varphi(\tilde\gamma(t))|
+\|\psi(\gamma(t))
    - \psi(\tilde\gamma(t))\|)\,d\mu(t)>\varepsilon^2/8(1+4\Lip(f))
$$
whenever $\gamma\in\Gamma(X)$ and
$\|\gamma-\hat\gamma\|_{\le k+1}<\delta$.
\endproclaim 

\demo{Proof} By Lemma~\ref{L1.1} we can
find $k\ge l$, $\delta_0>0$
and a set $A\subset T$ with $\mu(A)<\varepsilon/4$
such that for every $t\in T\setminus A$ and $0< r <\delta_0$
\begin{eqnarray}
\int_{Q_k(t,r)}
|\varphi(\tilde\gamma(s))-\varphi(\tilde\gamma(t))|\,d\mu(s)
&<&\varepsilon\mu(Q_k(t,r))/16\Lip(f),
\label{p1-lem11}\\
\noalign{\noindent\textrm{and}}
\int_{Q_k(t,r)}
\|\psi(\tilde\gamma(s))-\psi(\tilde\gamma(t))\|\,d\mu(s)
&<&\varepsilon\mu(Q_k(t,r))/4.
\label{p2-lem11}
\end{eqnarray}

For every $t\in S\setminus A$ and $n\in\N$ we use
Lemma~\ref{nonfpt-lem}
with 
$$\eta_n =\min(\eta,\delta_0)/(n+1)$$
to find $0<r_{t,n}<\eta_n$, $\delta_{t,n}>0$
and $\hat\gamma_{t,n}\in\Gamma_{k+1}(X)$
such that
$\|\hat\gamma_{t,n}-\tilde\gamma\|_{\le k}<\eta_n$,
$\hat\gamma_{t,n}(s)=\tilde\gamma(s)$ for
$s\in T\setminus Q_k(t,r_{t,n})$  and
either
\begin{equation}\label{q1-lem11}
\mu(Q_k(t,r_{t,n})\setminus\gamma^{-1}(E))
\ge\varepsilon\mu(Q_k(t,r_{t,n}))/8\Lip(f)
\end{equation}
or
\begin{equation}\label{q2-lem11}
\int_{Q_k(t,r_{t,n})\cap\gamma^{-1}(E)}
\|D_f(\gamma(s)) -D_f(\tilde\gamma(t))\|\,
  d\mu(s)\ge\varepsilon\mu(Q_k(t,r_{t,n}))/2,\qquad
\end{equation}
whenever $\gamma\in\Gamma(X)$ and
$\|\gamma(s)-\hat\gamma_{t,n}(s)\|
+\|D_{k+1}(\gamma(s)-\hat\gamma_{t,n}(s))\| <\delta_{t,n}$
for all
$s\in Q_k(t,r_{t,n})$.

The reason to introduce here the extra parameter $n$ is to enable us
to use the Vitali covering theorem later on
(see e.g.\ the definition of $\eta_n$).

If (\ref{q1-lem11}) holds, it follows using (\ref{p1-lem11})
and the fact that $\varphi(\tilde\gamma(t))=1$ and that
$\varphi$ vanishes outside $\gamma^{-1}(E)$ \pagebreak that
\begin{eqnarray*}
\lefteqn{\int_{Q_k(t,r_{t,n})}
|\varphi(\gamma(s))-\varphi(\tilde\gamma(s))|\,d\mu(s)}\\
&\ge&\int_{Q_k(t,r_{t,n})}
|\varphi(\gamma(s))-\varphi(\tilde\gamma(t))|\,d\mu(s)
-\int_{Q_k(t,r_{t,n})}
|\varphi(\tilde\gamma(s))-\varphi(\tilde\gamma(t))|\,d\mu(s)
\\
&>& \varepsilon\mu(Q_k(t,r_{t,n}))/8\Lip(f)
 -\varepsilon\mu(Q_k(t,r_{t,n}))/16\Lip(f)\\
&=&\varepsilon\mu(Q_k(t,r_{t,n}))/16\Lip(f).
\end{eqnarray*}
If (\ref{q2-lem11}) holds, we get using (\ref{p2-lem11}) that
$$
\int_{Q_k(t,r_{t,n})}
\|\psi(\gamma(s))
    - \psi(\tilde\gamma(s))\|\,d\mu(s)
 > \varepsilon\mu(Q_k(t,r_{t,n}))/4.
$$
Hence, in any case,
\begin{eqnarray}\label{p5-lem11}
\lefteqn{\int_{Q_k(t,r_{t,n})}
(|\varphi(\gamma(s))-\varphi(\tilde\gamma(s))|
+\|\psi(\gamma(s))
    - \psi(\tilde\gamma(s))\|)\,d\mu(s)}
\phantom{\varepsilon/4(1+4\Lip(f))}\\
&&>\varepsilon\mu(Q_k(t,r_{t,n}))/4(1+4\Lip(f))
\phantom{xxxxxx}\nonumber
\end{eqnarray}
for all $t\in S\setminus A$, $n\in\N$ and
$\gamma\in\Gamma(X)$ such that
$$\|\gamma(s)-\hat\gamma_{t,n}(s)\|
+\|D_{k+1}(\gamma(s)-\hat\gamma_{t,n}(s))\| <\delta_{t,n}$$
for all
$s\in Q_k(t,r_{t,n})$.

By the Vitali covering theorem, there are
$t_1,\dots,t_j$ and $n_1,\dots,n_j$ such that
$\{Q_k(t_i,r_{t_i,n_i})\}_{1\le i\le j}$ are disjoint and cover the set
$S\setminus A$
up to a set of measure less than $\varepsilon/4$.
Define now $\hat\gamma(s)=\hat\gamma_{t_i,n_i}(s)$
if $s\in Q_k(t_i,r_{t_i,n_i})$, $1\le i\le j$, and
$\hat\gamma(s)=\tilde\gamma(s)$ otherwise.
Clearly $\|\hat\gamma-\tilde\gamma\|_{\le k}<\eta$.
Letting $\delta=\min_{1\le i\le j}\delta_{t_i,n_i}$
then for every $\gamma\in\Gamma(X)$ with
$\|\gamma-\hat\gamma\|_{\le k+1}<\delta$
we get by adding~(\ref{p5-lem11}) over all the $j$ cubes
$Q_k(t_i,r_{t_i,n_i})$ that
\begin{eqnarray*}
\lefteqn{\int_T (|\varphi(\gamma(s))-\varphi(\tilde\gamma(s))|
+\|\psi(\gamma(s))
    - \psi(\tilde\gamma(s))\|)\,d\mu(s)}\\
&&>\varepsilon(\mu(S)-\varepsilon/2)/4(1+4\Lip(f))
\ge \varepsilon^2/8(1+4\Lip(f)).\\
\noalign{\vskip-36pt}
\end{eqnarray*}
\enddemo

We are now ready to prove the main theorem.

\proclaim{Theorem}\label{dif-gen2}
Suppose that $G$ is an open subset of a separable Banach space
$X${\rm ,} $\ls$ a norm separable subspace of $\Lin(X,Y)${\rm ,} and
$f:G\to Y$ is a Lipschitz function.
Then $f$ is Fr{\rm \'{\it e}}chet differentiable at $\Gamma$\/{\rm -}\/almost every
point $x\in X$ at which it is regular{\rm ,}
G{\rm \^{\it a}}teaux differentiable and $D_f(x)\in\ls$.
\endproclaim 

\demo{Proof} We may assume that $Y$ is separable and continue to use the
notation of the previous lemmas.
By Lemma~\ref{meas3-lem} there is
a residual subset $H$ of $\Gamma(X)$ such that the
restrictions of
$\Phi$ and $\Psi$ to $H$ are continuous.

Fix $\varepsilon>0$ and put
$$N =\{x\in E:
\limsup_{u\to 0} \|f(x+u)-f(x)-D_f(x)u\|/\|u\|>\varepsilon\}.$$
To prove the theorem it suffices to show that the set
$$M=\{\gamma\in H: \mu\{t:\gamma(t)\in N\}>2\varepsilon\}$$
is nowhere dense in $\Gamma(X)$. Assume
that this is not the case, then we can find
a nonempty open set $U$ in the closure of $M$.
Let
$\gamma_0\in M\cap U$
and making $U$ smaller if necessary we can
achieve that
\begin{equation}\label{ctr1}
\|\Phi(\gamma)-\Phi(\gamma_0)\|_{L_1}
+\|\Psi(\gamma)-\Psi(\gamma_0)\|_{L_1}
<\varepsilon^2/16(1+4\Lip(f))
\end{equation}
for every $\gamma\in U\cap H$.

Let $l_0\in\N$ and $\eta_0>0$ be such that
every $\gamma$ with
$\|\gamma-\gamma_0\|_{\le l_0}<3\eta_0$ belongs to $U$.
By Lemma~\ref{L1.2} we can find
$l\ge l_0$ so that
\begin{itemize}
\item
$\|\gamma_0^{l,s}-\gamma_0\|_{\le l_0}<\eta_0$ for every $s\in T$,
\item 
$\mu\{s\in T:
\|\Phi(\gamma_0^{l,s})-\Phi(\gamma_0)\|_{L_1}
\ge \varepsilon^2/32(1+\Lip(f))\}<\varepsilon/2$, and
\item $\mu\{s\in T:
\|\Psi(\gamma_0^{l,s})-\Psi(\gamma_0)\|_{L_1}
\ge \varepsilon^2/32(1+\Lip(f))\}<\varepsilon/2$.
\end{itemize}

By Fubini's theorem, there is $s\in T$ such that
$\tilde\gamma=\gamma_0^{l,s}$ has the properties that
$\mu\{t:\tilde\gamma(t)\in N\}>\varepsilon$
and
\begin{equation}\label{ctr1a}
\|\Phi(\tilde\gamma)-\Phi(\gamma_0)\|_{L_1}
+\|\Psi(\tilde\gamma)-\Psi(\gamma_0)\|_{L_1}
<  \varepsilon^2/16(1+\Lip(f)).
\end{equation}

By Lemma~\ref{nonf-lem},
we can find $k\ge l$, $\delta>0$ and
$\hat\gamma\in\Gamma_{k+1}(X)$
such that $\|\hat\gamma-\tilde\gamma\|_{\le k}<\eta_0$
and
\begin{equation}\label{ctr2}
\int (|\varphi(\gamma(t))-\varphi(\tilde\gamma(t))|
+\|\psi(\gamma(t))
    - \psi(\tilde\gamma(t))\|)\,d\mu(t)>\varepsilon^2/8(1+4\Lip(f))
\end{equation}
whenever $\gamma\in\Gamma(X)$ and
$\|\gamma-\hat\gamma\|_{\le k+1}<\delta$.
There are $\gamma\in U\cap H$ with\break
$\|\gamma-\hat\gamma\|_{\le k+1}<\delta$, and for any such
$\gamma$ both (\ref{ctr1a}) and (\ref{ctr2}) hold; hence
$$\|\Phi(\gamma)-\Phi(\gamma_0)\|_{L_1}
+\|\Psi(\gamma)-\Psi(\gamma_0)\|_{L_1}
>\varepsilon^2/16(1+4\Lip(f)).
$$
This contradicts (\ref{ctr1}) and thus the assumption that the
closure of $M$ contains a nonempty open subset led to a contradiction.
\enddemo

\proclaim{{C}orollary}
Assume that $X^\star$ is separable. Then any convex continuous
function on an open subset of $X$ is $\Gamma${\rm -}\/almost everywhere
Fr{\rm \'{\it e}}chet differentiable.
\endproclaim 

\demo{Proof} This is an immediate consequence of
Theorem~\ref{gat-ae}, Proposition~\ref{conv-reg} and Theorem~\ref{dif-gen2}.
\enddemo 

\proclaim{{C}orollary}\label{aspl-cor}
Assume that $X^\star$ is separable. Then any Lipschitz
function $f$ from an open subset $G$ of $X$ into $\R$
is $\Gamma$\/{\rm -}\/almost everywhere
Fr{\rm \'{\it e}}chet differentiable if and only if every $\sigma$ porous
set in $X$ is $\Gamma$\/{\rm -}\/null.
\endproclaim 

\demo{Proof} The ``if'' part is an immediate consequence of
Theorem~\ref{gat-ae}, Proposition~\ref{por-prp} and Theorem~\ref{dif-gen2}.

The ``only if'' part is trivial: As already remarked in
Section~\ref{sec:I}, if $A$ is porous, the Lipschitz function
$f(x)=\dist(x,A)$ is nowhere Fr\'echet differentiable on $A$.
\enddemo

\section{Spaces in which $\sigma$-porous sets are $\Gamma$-null}
\label{sec:S}

In this section we present the second main result of this paper
in which we identify a class of Banach spaces in which
$\sigma$-porous sets are $\Gamma$-null. The basic observation
behind these results is that,
if a surface passes through a point of a porous set $E$, then its
suitable small
deformation passes
through the center of a relatively large
ball that completely avoids $E$. This process may be iterated in order
to construct surfaces avoiding more and more of $E$. Unfortunately,
the final deformation cannot be guaranteed to remain small.
As we shall see in Corollary~\ref{ell-p}, in $\ell_2$ this problem is
unsurmountable. However, in $c_0$ any combination of small
perturbations is still small provided that different perturbations
use disjoint sets of coordinates, and the idea can be used to
show that sets porous
in the direction of all subspaces $\{x\in c_0: x_1=\dots=x_k=0\}$
are $\Gamma$-null. A simple decomposition of porous sets
(Lemma~\ref{pz-lem}) will then finish the proof
that every porous subset of $c_0$
is $\Gamma$-null. Moreover,
a close inspection of the first argument reveals
that the structure of $c_0$ is needed only asymptotically
(in the sense of the following definition), which will enable us to
extend the above arguments to several other spaces.

\numbereddemo{Definition}  \label{def:as}
Let $X$ be a Banach space and
$\{X_k\}_{k=1}^\infty$ a decreasing sequence of subspaces
of $X$.
The sequence $X_k$ of subspaces is said to be
asymptotically $c_0$ if
there is $C<\infty$ so that for every $n\in\N$,
$$
\begin{array}{cc}
(\exists k_1\in\N)(\forall u_1\in X_{k_1})
(\exists k_2\in\N)(\forall u_2\in X_{k_2})\dots
(\exists k_n\in\N)(\forall u_n\in X_{k_n})\\[5pt]
\|u_1+\dots+u_n\|\le C\max(\|u_1\|,\dots,\|u_n\|).
\end{array}
$$
\enddemo

\advance\eqcount by 14
A formally stronger requirement than that of Definition~\ref{def:as}
which is however easily seen to be equivalent to it is
that there is $C<\infty$ such
that for every $n\in\N$,
\begin{equation}\label{asc1}
\begin{array}{cc}
(\exists k_1\in\N)(\forall U_1 \subspace X_{k_1})
\dots
(\exists k_n\in\N)(\forall U_n \subspace X_{k_n})\\
(\forall u_i\in U_i)
\|u_1+\dots+u_n\|\le C\max(\|u_1\|,\dots,\|u_n\|),
\end{array}
\end{equation}
where the symbol $\subspace$ means
``a finite-dimensional subspace of''.

The main part of the proofs of the results of this section is contained
in the following lemma.

\proclaim{Lemma}\label{asc-lem}
Suppose that a sequence $\{X_k\}_{k=1}^\infty$ of subspaces of $X$
is asymptotically $c_0$. Then for every $c>0$
every set $E\subset X$ which is
$c$\/{\rm -}\/porous in the direction of all the subspaces $X_k$ is
$\Gamma$\/{\rm -}\/null.
\endproclaim

The notion of a set being
$c$-porous in the direction of a subspace is defined
in Section~\ref{sec:I}.   

\demo{Proof} It suffices to find a contradiction from the assumption
that there is a nonempty open subset $H$ of $\Gamma(X)$
having the property that, for some $\varepsilon>0$,
every nonempty open subset $G$ of $H$ contains a $\gamma\in\Gamma(X)$
such that
\begin{equation}\label{pp.assume}
\mu\{t: \gamma(t)\in E\} > \varepsilon.
\end{equation}

Let $\tilde\gamma_1\in H$ be such that
$\mu\{t:\tilde\gamma_1(t)\in E\}> \varepsilon$.
Find $\tilde m\in\N$ and $\delta_1>0$ such that
every $\gamma\in\Gamma(X)$ satisfying
$\|\gamma-\tilde\gamma_1\|_{\le \tilde m}<2\delta_1$
belongs to $H$.
We fix $m\ge\tilde m$ such that
$\|\tilde\gamma_1^{m,t}-\tilde\gamma_1\|_{\le\tilde m}<\delta_1$
for all $t\in T$ and use Fubini's theorem to find
$\bar t\in T$ such that the surface
$\gamma_1=\tilde\gamma_1^{m,\bar t}$
satisfies
$\mu\{t:\gamma_1(t)\in E\}> \varepsilon$.
We choose $M<\infty$ with the property
that for every $\gamma\in\Gamma_m(X)$ such that
$\|\gamma-\gamma_1\|_{\le m}\le\delta_1$,
\begin{equation}\label{pp.5}
\|\gamma(t)-\gamma(s)\|\le
M\max_{1\le j\le m}|t_j-s_j|
\textrm{ for all $t,s\in T$.}
\end{equation}

Denote $\kappa=c/4M$ and $K=4\max(\kappa,C/\delta_1)$,
where $C$ is the constant for which~(\ref{asc1}) holds
and denote
$$Q(s,r)=\{t\in\ell_\infty: |t_j-s_j|\le r\textrm{ for }
 j=1,\dots,m\}.$$
Choose $n\in\N$ so that
\begin{equation}\label{j.1}
(n-1)\varepsilon (\kappa/K)^m>2.
\end{equation}
It is this $n$ for which we intend to use~(\ref{asc1}).

For $i=1,\dots,n$ we will define inductively
indices $k_i\in\N$, surfaces
$\gamma_i,\psi_i\in\Gamma_m(X)$, finite sets
$S_i\subset T$, finite-dimensional subspaces
$U_i$ of $X_{k_i}$, sets $W_i,Q_i\subset T$, and
numbers $\delta_i>0$ so that in particular
the following statements hold
for $1\le i\le n$:
\begin{eqnarray}
&&\mu(W_i) >\varepsilon \label{pp.0a}\\
&&B(\gamma_i(s)+\psi_i(s),c\|\psi_i(s)\|)\cap E=\emptyset
\textrm{ if $s\in S_i$}\phantom{xxxx}\label{pp.3}
\end{eqnarray}
and the following statements hold
for $1\le i\le n-1$:
\begin{eqnarray}
&&\delta_{i+1}=\textstyle \frac{1}{4}
\min_{1\le l\le i}(\delta_l,
  c\min_{s\in S_l}\|\psi_l(s)\|) ,\phantom{xxxxxxx}\label{pp.0}\\
&&\|\gamma_{i+1}-\gamma_i\|_0<\delta_i, \label{pp.1a}\\
&&\|\gamma_{i+1}-(\gamma_i+\psi_i)\|_{\le m}<\delta_{i+1},
  \label{pp.1}\\
&&\|\gamma_{i+1}-\gamma_1\|_{\le m}<\delta_1. \label{pp.4}
\end{eqnarray}

Assume that for some $1\le i\le n$ we have already defined
$\gamma_i$, $\delta_i$ and all $\gamma_j$, $\delta_j$
$k_j$, $\psi_j$, $W_j$, $Q_j$, $S_j$, $U_j$ for $1\le j <i$.
Since $\gamma_1$ and $\delta_1$ have been already defined,
this is certainly true for $i=1$. We
show how to choose $k_i$, $W_i$, $Q_i$, $S_i$, $U_i$,
$\psi_i$ and, provided that $i<n$, $\delta_{i+1}$ and $\gamma_{i+1}$.

We find $k_i\in\N$ as in~(\ref{asc1}).
Let $W_i=\{t\in T:\gamma_i(t)\in E\}$ and
$$Q_i=\bigcup_{l=1}^{i-1} \bigcup_{s\in S_l}
 Q(s,\kappa\|\psi_l(s)\|).$$
The $\gamma_i$ has been chosen from among the surfaces for
which~(\ref{pp.assume}) holds, hence $\mu(W_i) > \varepsilon$ as
required by~(\ref{pp.0a}).
We show next that $W_i\cap Q_i=\emptyset$.
Since $Q_1=\emptyset$, there is nothing to prove for $i=1$.
Let $i>1$ and $t\in Q(s,\kappa\|\psi_l(s)\|)$ for some
$1\le l <i$.
We deduce from~(\ref{pp.4}) that~(\ref{pp.5}) is applicable to
$\gamma_i$ and infer, using also~(\ref{pp.1a})
and~(\ref{pp.1}), that
\begin{eqnarray*}
\lefteqn{\|\gamma_i(t)-(\gamma_l(s)+\psi_l(s))\|}\phantom{xx}\\
&\le&
\|\gamma_i(t)-\gamma_i(s)\|
+\|\gamma_i(s)-\gamma_{l+1}(s)\|
+\|\gamma_{l+1}(s)-(\gamma_l(s)+\psi_l(s))\|\\
&\le& M\kappa\|\psi_l(s)\|
+\sum_{j=l+1}^{i-1} \delta_j +\delta_{l+1}\\
&< & c \|\psi_l(s)\|.
\end{eqnarray*}
Hence~(\ref{pp.3}) implies that $\gamma_i(t)\notin E$
and thus $t\notin W_i$.

For every $t\in W_i$ choose
$u_{j,i}(t)\in X_{k_i}$ such that
$0<\|u_{j,i}(t)\|<\delta_i/4^j$ and
$B(\gamma_i(t)+u_{j,i}(t), c\|u_{j,i}(t)\|)\cap E=\emptyset$.
This can be done since, by assumption, $E$ is $c$-porous in the
direction of $X_{k_i}$.
Using Vitali's covering theorem (in $[0,1]^m$),
we find a finite set
$S_i\subset W_i$ and for each $s\in S_i$
vectors $u_i(s)=u_{j_i,i}(s)$ such that
the cubes $Q(s,K\|u_i(s)\|)$  are disjoint,
contained in $T\setminus Q_i$ and so that
\begin{equation}\label{pp.add}
\mu(\bigcup_{s\in S_i} Q(s,K\|u_i(s)\|)) \ge \mu(W_i)/2.
\end{equation}

We define next $U_i$ to be the span of $u_i(s)$, $s\in S_i$.
For each $s\in S_i$ we choose
$\omega_s\in\Gamma_m(\R)$ depending on the first $m$ coordinates
(so essentially $\omega_s:[0,1]^m\to\R$) such that
$0\le\omega_s\le 1$ and
\begin{eqnarray*}
\omega_s(t)&=&0\textrm{ for }t\notin Q(s,K\|u_i(s)\|),\\
\omega_s(s)&=&1\textrm{, and}\\
\|D_j\omega_s\|_0&\le& 2/(K\|u_i(s)\|)\textrm{ for }j=1,\dots,m.
\end{eqnarray*}

Let $\psi_i(t)$ be defined by
$\psi_i(t)=\sum_{s\in S_i} \omega_s(t)u_i(s)$.
For each $t$, $\psi_i(t)\in U_i$ and the same is true for the
partial derivatives of $\psi_i$. Since $\psi_i(s)=u_i(s)$
for $s\in S_i$, (\ref{pp.3}) follows from the choice of the $u_i(s)$.
The definition of $u_i(s)$ gives also that
$\|\psi_i\|_0 < \delta_i/4$ and consequently,
 from (\ref{pp.1}) and \pagebreak $\delta_{l+1}\le\delta_l/4$,
\begin{eqnarray*}
\|\gamma_i+\psi_i -\gamma_1\|_0 &\le&
\|\sum_{l=1}^i \psi_l\|_0 +\sum_{l=1}^{i-1}
\|\gamma_{l+1}-(\gamma_l+\psi_l)\|_0\\ &\le&
2\sum_{l=1}^{i} \delta_l/4 < \delta_1.
\end{eqnarray*}
Since $D_j\psi_l(t)\in U_l$ for all $1\le j \le m$ and $t\in T$,
we get from (\ref{asc1}) that
$$\| D_j(\sum_{l=1}^i \psi_l)(t)\|
\le C\max_{1\le l\le i} \|D_j\psi_l(t)\|
\le 2C/K.$$
Hence, by (\ref{pp.1}),
\begin{eqnarray}\label{pp.l1}
\|\gamma_i+\psi_i -\gamma_1\|_j
&\le&
\|\sum_{l=1}^i \psi_l\|_j
+\sum_{l=1}^{i-1} \|\gamma_{l+1}-(\gamma_l+\psi_l)\|_j\\
&\le& 2C/K + \sum_{l=1}^{i-1} \delta_{l+1} < \delta_1,\nonumber
\end{eqnarray}
and thus
$\gamma_i+\psi_i\in H$.

If $i=n$, this finishes the construction. If $i<n$, we
still have to define $\delta_{i+1}$ and  $\gamma_{i+1}$ and show
that~(\ref{pp.0}), (\ref{pp.1a}), (\ref{pp.1}) and~(\ref{pp.4})
remain valid. 
Since $\psi_i(s)\ne 0$, (\ref{pp.0}) may be used as a definition of
$\delta_{i+1}$.
By our assumption on $H$, there is a
$\tilde\gamma_{i+1}\in H$ such that
$$\|\tilde\gamma_{i+1} -(\gamma_i+\psi_i)\|_{\le m}
< \delta_{i+1}
\textrm{ and }
\mu(\{t: \tilde\gamma_{i+1}(t)\in E\})>\varepsilon.$$
By Fubini's theorem there is a $\bar t\in T$ such that the surface
$$\gamma_{i+1}(t)=
\tilde\gamma_{i+1}^{m,\bar t}$$
satisfies
$\mu\{t:\gamma_{i+1}(t)\in E\}> \varepsilon$.
Since $\gamma_i$ and $\psi_i$ belong to $\Gamma_m(X)$, it
follows that~(\ref{pp.1a}) holds. Since $\|\psi_i\|_0 < \delta_i/4$,
also~(\ref{pp.1}) holds.
To establish~(\ref{pp.4}), we use~(\ref{pp.l1}) to get
\begin{eqnarray*}
\|\gamma_{i+1} -\gamma_1\|_{\le m}
&\le& \|\gamma_{i+1}-(\gamma_i+\psi_i)\|_{\le m}
  + \|\gamma_i+\psi_i -\gamma_1\|_{\le m}\\
&\le& \delta_{i+1} + 2C/K + \sum_{l=1}^{i-1} \delta_{l+1}
< \delta_1.
\end{eqnarray*}
This finishes the $i$-th step of the construction.
We now loop back to the choice of $k_{i+1}$ according to~(\ref{asc1})
and continue the inductive process.

Finally, by~(\ref{pp.0a}) and~(\ref{pp.add}) we get
\begin{eqnarray*}
\mu(Q_{i+1})
&=&\mu(Q_i)+
 \sum_{s\in S_i} \mu(Q(s,\kappa\|u_i(s)\|))\\
&=&\mu(Q_i)+
 (\kappa/K)^m \sum_{s\in S_i} \mu(Q(s,K\|u_i(s)\|))\\
&\ge & \mu(Q_i) + (\kappa/K)^m \varepsilon/2.
\end{eqnarray*}
Hence $\mu(Q_{n})\ge (n-1)(\kappa/K)^m \varepsilon/2 > 1$
by our choice of $n$ in (\ref{j.1})
and we get the desired contradiction.
\enddemo

In order to apply Lemma~\ref{asc-lem}
we need the following lemma which is a variant of Lemma~4.6
from~\cite{pz2}.

\proclaim{Lemma}\label{pz-lem}
Let $U,V$ be subspaces of a Banach space $X$
such that for some $\eta<\infty$ every
$x\in U+V$ may be written as $x=u+v$ where
$u\in U${\rm ,} $v\in V$ and $\max(\|u\|,\|v\|)\le\eta \|x\|$.
Then{\rm ,} if a set $E \subset X$ is $c$\/{\rm -}\/porous in
the direction
$U+V${\rm ,} then we can write $E = A \cup B${\rm ,}
where
$A$ is $\sigma$\/{\rm -}\/porous in direction $U$ and $B$ is
$c/2\eta$\/{\rm -}\/porous in direction $V$.
\endproclaim

\demo{Proof} Denote by $B_m$ the set of those $x\in E$
for which there is a $v\in V$ with
$\|v\|< 1/m$ and $B(x+v, c\|v\|/2\eta)\cap E=\emptyset$.
Clearly, $B=\bigcap_{m=1}^\infty B_m$ is $c/2\eta$-porous in
direction $V$. Thus it is
sufficient to prove that
each set $E\setminus B_m$ is porous
in direction $U$. 

Let $x \in E\setminus B_m$ and
$0 < \varepsilon < 1/m$.
Since $E$ is $c$-porous in
direction $U+V$, we can find $z \in U+V$ such that
$0<\|z\|<\varepsilon/\eta$ and $B(x+z,c\|z\|) \cap E =
\emptyset$. Write $z=u+v$ with $u\in U$, $v\in V$
and $\max(\|u\|,\|v\|)\le\eta \|z\|$.

It suffices to show that
$B(x+u,c\|u\|/2\eta) \cap (E\setminus B_m) =
\emptyset$.
For this assume that
$y \in B(x+u,c\|u\|/2\eta) \cap (E\setminus B_m)$,
then
$B(y+v,c\|v\|/2\eta)
\subset B(x+u+v,c\|u+v\|) = B(x+z,c\|z\|)
\subset X\setminus E$.
Since $\|v\| \le \eta\|z\| < \frac{1}{m}$
this shows that $y \in B_m$ which contradicts
our assumption.
\enddemo

We can now give concrete examples of spaces having the property that
every $\sigma$-porous subset in them is $\Gamma$-null.

\proclaim{Proposition} \label{j1-prop}
Assume that $X$ has a Schauder basis $\{x_i\}_{i=1}^\infty$
such that for a suitable sequence $\{n_k\}_{k=1}^\infty$
the sequence $X_k=\clspan \{x_i\}_{i > n_k}$ is asymptotically $c_0$.
Then any $\sigma$\/{\rm -}\/porous set $E$ in $X$ is $\Gamma$\/{\rm -}\/null.
\endproclaim 

\demo{Proof} Let $U_k=\lspan \{x_i\}_{i \le n_k}$. Then clearly $X=U_k\oplus X_k$
for every $k$. If $E$ is porous we can write (by Lemma~\ref{pz-lem})
for every $k$, $E=A_k\cup B_k$ where $A_k$ is $\sigma$-porous
in the direction $U_k$ and $B_k$ is $c_1$-porous
in the direction $X_k$ for some fixed $0<c_1<1$. By
Lemma~\ref{asc-lem}, $\bigcap_{k=1}^\infty B_k$ is $\Gamma$-null.
It remains to show that every porous (and thus $\sigma$-porous) set
$A$ in direction of a finite-dimensional subspace is $\Gamma$-null.
A simple compactness argument shows that such a set $A$ is actually
directionally porous. Hence the Lipschitz function from $X$ to $\R$
defined by $f(x) = \dist(x,A)$ is not G\^ateaux differentiable at any
point of $A$. By Theorem~\ref{gat-ae} $A$ is $\Gamma$-null.
\enddemo

We recall now (see e.g.\ \cite[Lemma 2.13]{jlps}) the simple fact that
if $W$ is a finite co-dimensional subspace of a Banach space $X$, then
there is a finite-dimensional subspace $V$ of $X$ so that every $x\in
X$ can be written as $x=v+w$ with $v\in V$, $w\in W$ and
$\|v\|\le 2\|x\|$, $\|w\|\le 3\|x\|$. By using this fact the
argument of Proposition~\ref{j1-prop} can be carried through for any space
admitting an asymptotically $c_0$ sequence of finite co-dimensional
subspaces. In particular,
any subspace of $c_0$ also has the property that any
porous set in it is $\Gamma$-null.

The next lemma allows us to do some iteration arguments concerning the
property we are interested in here.

\proclaim{Lemma}\label{st-prop}
Consider the following property of a Banach space $X$.
\begin{itemize}
\item[$(\star)$]
Whenever $X$ is a complemented subspace of $Y$ and $E\subset
Y$  is $\sigma$\/{\rm -}\/porous in the direction $X${\rm ,}
then $E$ is $\Gamma$\/{\rm -}\/null.
\end{itemize}

Any finite or $c_0$ {\rm (}\/infinite\/{\rm )} direct sum of spaces having $(\star)$
also has $(\star)$.

Moreover{\rm ,} if there are
subspaces $U_k,V_k$ of $X$ such that
$U_k$ are complemented in $X$
and have property $(\star)$ while the sequence $V_k$
is asymptotically $c_0$ and there is $\eta<\infty$
such that every $x\in X$ can be written as
$x=u+v$ where $u\in U_k${\rm ,} $v\in V_k$
and $\max(\|u\|,\|v\|)\le\eta\|x\|${\rm ,} then
$X$ has property $(\star)$.
\endproclaim

\demo{Proof} Assume that $X_1$ and $X_2$ have $(\star)$,
$X=X_1\oplus X_2$ is complemented in $Y$ and $E\subset Y$ is
$\sigma$-porous in the direction $X$.
By
Lemma~\ref{pz-lem}, $E=E_1\cup E_2$, where $E_i$ is
$\sigma$-porous in the direction $X_i$. It follows from
$(\star)$ that each $E_i$ is $\Gamma$-null and thus $E$ is
$\Gamma$-null. Hence $(\star)$ is stable under finite direct
sums.

We prove next the ``moreover'' statement. This statement and the
stability of $(\star)$ under finite direct sums imply that
$(\star)$ is also closed under $c_0$ direct sums.

Let $E$ be a $c$-porous subset of $Y$ in the direction $X$.
By Lemma~\ref{pz-lem},
$E= A_k\cup B_k$, where $A_k$ is $\sigma$-porous in the
direction $U_k$ and $B_k$ is
$c/2\eta$-porous in the direction $V_k$. By
assumption the $A_k$ are $\Gamma$-null.
Since the sequence $V_k$ is asymptotically $c_0$,
and the set $B=\bigcap_{k=1}^\infty B_k$ is $c/2\eta$-porous
in the direction of every $V_k$, it follows from Lemma~\ref{asc-lem}
that $B$ is $\Gamma$-null and so is
$E=B\cup\bigcup_{k=1}^\infty A_k$. \phantom{myloom}
\enddemo

Recall that any compact countable Hausdorff topological space
$K$ is homeomorphic to the space $K_\alpha$ of ordinals $\le\alpha$
in the order topology for some countable ordinal $\alpha$. Recall also
that $C(K_\omega)$ is isomorphic to $c_0$ and that every $C(K_\alpha)$
is isomorphic to $(C(K_{\beta_1})\oplus\dots\oplus
C(K_{\beta_k})\dots)_{c_0}$ for suitable $\{\beta_k\}_{k=1}^\infty$
smaller than $\alpha$. Hence we get from Lemma~\ref{st-prop} by
(countable) transfinite induction that each such $C(K)$ has property
$(\star)$.

Summing up the preceding observations we get

\proclaim{Theorem}\label{c0-por}
The following spaces have the property that all their\break
$\sigma$\/{\rm -}\/porous subsets are $\Gamma$\/{\rm -}\/null\/{\rm :}
$c_0${\rm ,} $C(K)$ with $K$ compact countable{\rm ,} the
Tsirelson space{\rm ,} subspaces of $c_0$.
\endproclaim 

The Tsirelson space (as first defined in \cite{tsirelson}) is a space
with an unconditional basis which satisfies the assumption in
Proposition~\ref{j1-prop}. An important feature of this space (in general,
and for us here) is that it is reflexive.

 From Corollary~\ref{aspl-cor} we see that each of the spaces
listed in Theorem~\ref{c0-por} has the property
that every real-valued Lipschitz function defined on its open subset
is Fr\'echet differentiable $\Gamma$-almost everywhere. Stronger
results follow by application of Theorem~\ref{dif-gen2}. For optimal
results in this direction it is desirable to have
information when
$\Lin(X,Y)$ is separable for every space $Y$ with the RNP. Call
this property of a Banach space $(\star\star)$.

\proclaim{Lemma}\label{sl-lem}
The class of spaces $X$ with the property $(\star\star)$
is closed under finite direct sums and under infinite direct sums in
the sense of $c_0$.
\endproclaim

\demo{Proof} For finite direct sums this is obvious. Assume that
$X=(\bigoplus_{i=1}^\infty X_i)_{c_0}$ with each $X_i$ having
$(\star\star)$.
Denote by
$\pi_k$ the natural projection onto $\sum_{i=1}^{k^{\phantom{|}}} X_i$
and let $V_k=(I-\pi_k)X$, $k=1,2\dots$.
Let $L\in\Lin(X,Y)$ with $Y$ having the RNP, let $\varepsilon>0$
and let $0<\eta<1$ be such that $\eta/(1-\eta)<\varepsilon$.

Since $Y$ has the RNP, the set
$\{Lz: \|z\|\le 1\}\subset Y$ has slices with arbitrarily small diameter.
Thus there are $u^\star\in U^\star$ and a real $c$ so that the set
$$S=\{Lz: \|z\|\le 1, u^\star(Lz) > c\}$$
is nonempty and of diameter $<\eta$. Let $z\in Z$ with
$\|z\|\le 1$ be such that $Lz\in S$ and let $m$ be such that
$\|z-\pi_mz\|<\eta$. Whenever $w\in V_m$ with $\|w\|\le 1$ then
$\|z\pm(1-\eta)w\|\le 1$, so at least one of the vectors
$L(z\pm(1-\eta)w)$ also belongs to $S$. Hence
$\|L((1-\eta)w)\|<\eta$ or
$\|L-L\circ\pi_m\|<\eta/(1-\eta)<\varepsilon$. Since
$\Lin(\sum_{i=1}^m X_i,Y)$ is separable for every $m$, it follows
that $\Lin(X,Y)$ is separable.
\enddemo

Since $\R$ clearly has $(\star\star)$, it follows
from Lemma~\ref{sl-lem} that $c_0$ has
$(\star\star)$ and by countable transfinite induction that
$C(K)$ has $(\star\star)$ for every compact countable~$K$.

We show next that every subspace $Z$ of $c_0$
has property $(\star\star)$.
Let $V_k=\{x\in Z: x_1=\dots=x_k=0\}$.
By the observation from \cite{jlps}
used already above there is for every $k$ a finite-dimensional subspace
$U_k$ of $Z$ so that every $z\in Z$ can be written as $u+v$ with
$u\in U_k$, $v\in V_k$ and $\|u\|,\|v\|\le 3\|z\|$. For every
$L\in\Lin(Z,Y)$ and every $\varepsilon>0$ there is a $k$ such that
$\|Lv\|\le \varepsilon\|v\|$ for every $v\in V_k$. Indeed, if this
were false, we would get by a standard gliding bump argument that $L$
is an isomorphism on a subspace of $Z$ isomorphic to $c_0$ which
contradicts the assumption that $Y$ has the RNP.

Since $U_k$ is finite-dimensional, there is, given $\varepsilon>0$,
a countable set $\ls_k\subset\Lin(Z,Y)$ such that for every
$L\in\Lin(Z,Y)$ with $\|Lv\|\le\varepsilon\|v\|$ for $v\in V_k$
there is an $\hat L\in\ls_k$ such that
$\|\hat Lv\|\le\varepsilon\|v\|$ for $v\in V_k$
and $\|\hat Lu-Lu\|\le\varepsilon\|u\|$ for $u\in U_k$.
By decomposing each $z\in Z$ as above, we deduce that
$\|\hat Lz-Lz\|\le 6\varepsilon\|z\|$ for every such $z$.
Hence every $L\in\Lin(Z,Y)$ has distance ${}\le 6\varepsilon$
from the countable set $\bigcup_k\ls_k$ and this proves the
separability of $\Lin(Z,Y)$.

From Theorems~\ref{dif-gen2} and~\ref{c0-por} and the preceding
arguments we deduce

\proclaim{Theorem}\label{gam-dif} \hskip-8pt
The following spaces have the property that every Lipschitz
mapping of them into a space with the {\rm RNP} is Fr{\rm \'{\it e}}chet
differentiable $\Gamma$\/{\rm -}\/almost everywhere\/{\rm :}
$C(K)$ for compact countable $K${\rm ,} subspaces of $c_0$.
\endproclaim 

\demo{{R}emark} 
Theorem~\ref{gam-dif} does not hold for subspaces of $C(K)$,
$K$ countable. As remarked in \cite{jlps} the Schreier space
which is the completion of the space of eventually zero sequences
$\{a_i\}_{i=1}^\infty$ with respect to the norm
$$\|\{a_i\}\|=\sup\{\sum_{k=1}^n|a_{i_k}|: n\in\N
\textrm{ and } n\le i_1<i_2<\dots<i_n\}$$
is isomorphic to a subspace of $C(\omega^\omega)$ and there is a
Lipschitz map from this space into $\ell_2$ which is nowhere
Fr\'echet differentiable.
\enddemo

\section{The mean value theorem}\label{sec:M}

Let $X$ be a separable Banach space and $Y$ a space having the RNP.
Let $D$ be a convex open set in $X$ and $D_0$ be the subset
of $D$ consisting 
of points where the G\^ateaux derivative $D_f(x)$ exists. We know by
Theorem~\ref{Gat-orig} (resp.\ Theorem \ref{gat-ae})
that $D\setminus D_0$ is Gauss
null (resp.\ $\Gamma$-null).

Put for $u\in X$,
\begin{eqnarray*}
R_u&=&\{\frac{f(x+tu)-f(x)}{t}:x,x+u\in D
\textrm{ and } t>0\},\\
\tilde R_u &=& \{D_f(x)u: x\in D_0\}.
\end{eqnarray*}
Then a simple application of the separation theorem,
Theorem~\ref{Gat-orig}
and the property of Gauss null sets that whenever
$x,x+tu\in D$ then there is a point $x_0$ arbitrarily close to $x$ such
that $x_0+su\in D$ for almost all $0\le s\le t$, shows that
the closed convex hull of $\tilde R_u$ is equal to the closed convex
hull of $R_u$
(see e.g.\ \cite[Lemma 2.12]{jlps}).

The proofs in \cite{p} and \cite{lp} of existence of points of
Fr\'echet differentiability of Lipschitz functions $f:X\to\R$ where
$X^\star$ is separable actually show the following stronger
assertion: Let $D$ be a convex open set in $X$ and let $u,v\in D$
and a constant $m$ be such that
$f(v)-f(u) > m$. Then there is a point $x\in D$ in which $f$ is
Fr\'echet differentiable and $D_f(x)(v-u) > m$.

In view of the observation above this statement is equivalent to
saying that whenever $f$ is a real-valued Lipschitz map on an open set
$G$ in $X$ (with $X^\star$ separable) then any nonempty slice $S$ of
the set $\Upsilon$ of G\^ateaux derivatives of $f$ contains $D_f(x)$
where $x\in G$ and $f$ is Fr\'echet differentiable at $x$.
Recall that a slice $S$ of $\Upsilon$ is a set of the form
$S(\Upsilon,v,\delta)$ where $v\in X$, $\delta>0$ and
$$S(\Upsilon,v,\delta)
=\{T\in \Upsilon: Tv > \alpha -\delta,
  \alpha = \sup_{T\in \Upsilon} Tv\}.$$

A natural generalization of this assertion to maps into $\R^n$ or
more general spaces $Y$ would be as follows: Let $\Upsilon$ be again
the set of G\^ateaux derivatives of $f$ and consider any slice
$S=S(\Upsilon,\{v_i\}_{i=1}^m,\{y^\star_i\}_{i=1}^m,\delta)$
of $\Upsilon$ where $m\in\N$, $\{v_i\}_{i=1}^m\subset X$,
$\{y^\star_i\}_{i=1}^m\subset Y^\star$, $\delta>0$ and
$$S(\Upsilon,\{v_i\},\{y^\star_i\},\delta)
=\{T\in \Upsilon: \sum_{i=1}^m y_i^\star(T v_i) > \alpha -\delta,
  \alpha = \sup_{T\in \Upsilon} \sum_{i=1}^m y_i^\star(T v_i)\}.$$
Then $S$ contains a point of the form $D_f(x)$
where $f$ is Fr\'echet differentiable at~$x$.
Unfortunately, this generalization does not hold even for maps into
$\R^2$ as the following example shows.

\demo{Example {\rm 5.1 (\cite{pt})}} 
Let $1<p<\infty$ and let $n$ be an integer with $n>p$. Then there is a
Lipschitz map $f=(f_1,f_2,\dots,f_n)$ from $\ell_p$ into $\R^n$ such
that whenever $f$ is Fr\'echet differentiable at a point $x$ then
$\sum_{i=1}^n D_{f_i}(x)e_i=0$, where $\{e_j\}_{j=1}^\infty$ denote
the unit basis in $\ell_p$, but $f$ is G\^ateaux differentiable
at $x=0$ and we have that $\sum_{i=1}^n D_{f_i}(0)e_i=1$.
\enddemo
\advance\theoremcount by 1

The next theorem shows however, that in the sense of $\Gamma$-almost
everywhere the mean value theorem for Fr\'echet derivatives holds also
for maps into spaces of dimension greater than one.

\proclaim{Theorem}\label{div-est}
Suppose that $f:G\to Y$ is a Lipschitz mapping
which is Fr{\rm \'{\it e}}chet differentiable at\/ $\Gamma$\/{\rm -}\/almost
every point of an open subset $G$ of a Banach space
$X$. Then{\rm ,} for every slice $S$ of the set
of G{\rm \^{\it a}}teaux derivatives of $f${\rm ,} the set of points $x$ at which $f$ is
Fr{\rm \'{\it e}}chet differentiable and $D_f(x)\in S$ is not $\Gamma$\/{\rm -}\/null.
\endproclaim 

\demo{Proof} Let $S=S(\Upsilon,v_1,\dots,v_n,y_1^\star,\dots,y_n^\star,\delta)$,
where $\Upsilon$ is the set of G\^ateaux derivatives of $f$.
We can assume, without loss of generality, that
$\sum_{k=1}^n \|v_k\| = \sum_{k=1}^n \|y_k^\star\| = 1$.
Let $x_0\in G$ be such that $f$ is
G\^ateaux differentiable at $x_0$ and
$D_f(x_0)\in S(\Upsilon,v_1,\dots,v_n,y_1^\star,\dots,y_n^\star,\delta/4)$.
Define $f_0(y)=f(x_0)+D_f(x_0)(y-x_0)$ and find $r>0$ such that
$\|f(x_0+v)-f_0(x_0+v)\|\le \delta\|v\|/8$ for
$v\in\lspan\{v_k\}_{k=1}^n$ and $\|v\|\le r$.
Let $\gamma_0(t)=x_0+r\sum_{k=1}^n t_k v_k$ and
consider any $\gamma\in\Gamma(X)$ such that
$\|\gamma-\gamma_0\|_{\le n}<\delta r/8(1+\Lip(f))$.
For any $t\in T$ we have
\begin{eqnarray*}
\|f(\gamma(t))-f_0(\gamma_0(t))\|
&\le&
\|f(\gamma(t))-f(\gamma_0(t))\|+\|f(\gamma_0(t))-f_0(\gamma_0(t))\|\\
&\le&
\Lip(f)\|\gamma-\gamma_0\|_0+\delta\|\gamma_0(t)-x_0\|/8
 < \delta r/4.
\end{eqnarray*}
Since $\sum_{k=1}^n \|y_k^\star\| = 1$, integration with respect to the
first $n$ variables and the divergence
theorem imply that
$$ \int_T 
\sum_{k=1}^n \frac{\partial}{\partial t_k}
 \sprod{y_k^\star}{f(\gamma(t))-f_0(\gamma_0(t))}\, d\mu(t) <  \delta r/2.$$
It follows that for $t$ belonging  to a set of positive measure,
\begin{eqnarray*}
\sum_{k=1}^n \frac{\partial}{\partial t_k}
 \sprod{y_k^\star}{f(\gamma(t))}
&>& 
\sum_{k=1}^n \frac{\partial}{\partial t_k}
\sprod{y_k^\star}{f_0(\gamma_0(t))}  - r\delta /2 \\
&=&
r\sum_{k=1}^n \sprod{y_k^\star}{D_f(x_0)(v_k)}- r\delta /2.
\end{eqnarray*}
If, in addition, $\gamma$ is such that
$f$ is Fr\'echet differentiable
at $\gamma(t)$ for almost every~$t$, we conclude that
for $t$ belonging  to a set of positive measure,
\begin{eqnarray*}
\lefteqn{\sum_{k=1}^n \sprod{y_k^\star}{D_f(\gamma(t))(v_k)}
= \frac{1}{r}\sum_{k=1}^n
 \sprod{y_k^\star}{D_f(\gamma(t)) (D_k\gamma_0)}}\\
&\ge&
\frac{1}{r}
\sum_{k=1}^n \left(\frac{\partial}{\partial t_k}
 \sprod{y_k^\star}{f(\gamma(t))}
- \Lip(f)\|D_k\gamma(t)-D_k\gamma_0(t)\|\|y_k^\star\|\right)\\
&>&
\sum_{k=1}^n \sprod{y_k^\star}{D_f(x_0)(v_k)}-3\delta/4
\end{eqnarray*}
and thus $D_f(\gamma(t))\in S$ for $t$
from a set of positive measure.
\enddemo

\demo{{R}emark} 
The same proof shows that if we assume only that $f$ is G\^ateaux
differentiable $\Gamma$-almost everywhere and that we are given a
$\Gamma$-null set $N$ then every slice of the set of G\^ateaux
derivatives contains some $D_f(x)$
where $x\notin N$ and $f$ is Fr\'echet differentiable at $x$.
\enddemo

\proclaim{{C}orollary}\label{ell-p}
In $\ell_p$, $1<p<\infty${\rm ,} there are porous sets which are not
$\Gamma$\/{\rm -}\/null{\rm ,} and thus also real\/{\rm -}\/valued Lipschitz functions
whose sets of non Fr{\rm \'{\it e}}chet differentiability are not $\Gamma$\/{\rm -}\/null.
\endproclaim 

\demo{Proof} This is an immediate consequence of
Corollary~\ref{aspl-cor}, Example~5.1 and
Theorem~\ref{div-est}.
\enddemo

\demo{{R}emark} 
If $\pi$ is a projection of $X$ onto its subspace $Y$ then
a set $A\subset Y$ is $\Gamma$-null in $Y$ if and only if
$\pi^{-1}(A)$ is $\Gamma$-null in $X$. This follows by observing
that the map $\varphi:\Gamma(X)\to\Gamma(Y)$ defined by
$\varphi(\gamma)=\pi\circ\gamma$ is onto, hence open
(by the open mapping theorem), and hence it maps residual
sets onto residual sets. It follows that the statement of
Corollary~\ref{ell-p} holds also for any Banach
space which contains some $\ell_p$, $1<p<\infty$,
as a complemented subspace.
\enddemo

\section{Remarks and problems}\label{sec:R}

We start by stating the following differentiability conjecture.

\specialnumber{1}
\proclaim{Conjecture}\label{dc}
Assume that $X^\star$ is separable{\rm ,}
$\{g_i\}_{i=1}^\infty$ are
Lipschitz mappings from $X$ to $\R$
and $f_i$ are Lipschitz mappings from $X$ to $Y_i$
where each $Y_i$ has the {\rm RNP.} Then
there is a point $x\in X$ at which
all $g_i$ are Fr{\rm \'{\it e}}chet differentiable
and all $f_i$ are G{\rm \^{\it a}}teaux differentiable.
\endproclaim

It is unknown for which spaces $X$ this conjecture holds.
In particular, it is unknown if it holds for $X=\ell_2$.
It may hold even for every $X$ with $X^\star$ separable.
The results of Sections~\ref{sec:F} and \ref{sec:S}
show that it holds for $c_0$, $C(K)$ with $K$ compact countable, the
Tsirelson space and subspaces of $c_0$. These are the first
and so far only examples
of infinite-dimensional Banach spaces for which it is
known that the conjecture holds.

The differentiability conjecture is connected to another (well-known)
open problem. Assume that $X$ is Lipschitz equivalent to $Y$; i.e.\
there is a Lipschitz mapping $f$ from $X$ onto $Y$ so that $f^{-1}$
exists and is also Lipschitz. Must $X$ be linearly isomorphic to $Y$
if we assume that $X$ is separable and has the RNP?
It is trivial to check that if the Lipschitz equivalence is
Fr\'echet differentiable at some point $x_0$ then $D_f(x_0)$ is a
linear isomorphism from $X$ onto $Y$. If $f$ is only G\^ateaux
differentiable at $x_0$ then in general $D_f(x_0)$ will be only an
isomorphism from $X$ into $Y$. If however besides $f$ being G\^ateaux
differentiable at $x_0$ we also know that the functions
$g_i=y_i^\star\circ f$, $1\le i <\infty$ are Fr\'echet differentiable
at $x_0$, where $\{y_i^\star\}$ is a norming sequence in the unit ball
of $Y^\star$, then $D_f(x_0)$ must be surjective. All this is
discussed in \cite[Chap.~7]{BL}. Consequently,

\proclaim{Proposition} \label{iso-gen}
If $X$ has the {\rm RNP} and
satisfies the differentiability conjecture{\rm ,}
then every Banach space Lipschitz
equivalent to $X$ is linearly isomorphic
to $X$.
In particular{\rm ,} this holds if $X$ is the Tsirelson space.
\endproclaim 

In \cite[Chap.~7]{BL} it is proved for a large class of spaces $X$
that if $Y$ is Lipschitz
equivalent to $X$ then $Y$ is linearly isomorphic
to $X$. The Tsirelson space is the first known example where the
isomorphism can always be taken to be the G\^ateaux derivative of
the Lipschitz equivalence at a suitable point.

A similar situation occurs for Lipschitz quotient maps introduced and
studied in \cite{bjlps}.

\proclaim{Proposition} \label{quo-gen}
If $X$ is separable and reflexive and satisfies the
differentiability conjecture{\rm ,} then any Banach space $Y$
which is a Lipschitz
quotient of $X$ is already a linear quotient of $X$.
In particular{\rm ,} this holds for $X$ the Tsirelson space.
\endproclaim 

\demo{Proof} We prove that $Y$ is reflexive (and thus has the RNP).
Let $f$ be a  Lipschitz quotient map from $X$ onto $Y$,
It follows from \cite[Theorem~3.18]{bjlps} that $Y^\star$ is
separable. If $\{y_i^\star\}_{i=1}^\infty$ is a norm dense subset
of $Y^\star$, it follows from the
differentiability conjecture that that there is a point
$x_0\in X$ at which all the functions $y_i^\star\circ f$ are
Fr\'echet differentiable.  This easily implies that
$y^\star\circ f$ is Fr\'echet differentiable at $x_0$ for every
$y^\star\in Y^\star$.  The mapping which
assigns to every
$y^\star$ the Fr\'echet derivative of $y^\star\circ f$ at $x_0$ is
an isomorphic embedding (since $f$ is a Lipschitz quotient map).
Hence $Y^\star$ is isomorphic to a subspace of $X^\star$ and thus $Y$
is reflexive.
Using again the differentiability conjecture, we get that for some
$x_1\in X$ the G\^ateaux derivative of $f$ at $x_1$ exists and all
$y^\star\circ f$ $(y^\star\in Y^\star)$ are Fr\'echet differentiable at
$x_1$. An argument similar to that of\break \cite[Prop.~3.10]{bjlps}
will now show that $D_f(x_1)$ is a linear quotient map from $X$ onto~$Y$.
\enddemo

Another related problem is the following:

\numbereddemo{Problem}\label{pbm3}
Let $f$ be a Lipschitz equivalence from a separable Banach space $X$
with the RNP onto a Banach space $Y$. If $A\subset X$ is $\Gamma$-null,
does this imply that $f(A)$ is $\Gamma$-null?
\enddemo

For Haar null or Gauss null sets the answer to the analogue of this
question is negative even for $C^\infty$ diffeomorphisms
(see \cite[Chap.~5]{BL} for details and references).
For $\Gamma$-null sets the situation is considerably better;
since for any $C^1$ diffeomorphism of $X$ onto $Y$ the mapping
$\psi_f(\gamma) = f\circ\gamma$ is a homeomorphism of
$\Gamma(X)$ onto $\Gamma(Y)$, we have

\proclaim{Proposition} \label{prop:inv}
If $f$ is a $C^1$ diffeomorphism of a separable Banach space $X$
onto a separable Banach space $Y$ and $A\subset X$ is $\Gamma$\/{\rm -}\/null{\rm ,}
then $f(A)\subset Y$ is $\Gamma$\/{\rm -}\/null.
\endproclaim 

If the answer to Problem~\ref{pbm3} were positive it
would solve the problem whether Lipschitz equivalence implies
isomorphism (for separable spaces with the RNP). Indeed, this would
imply with the help of Theorem~\ref{gat-ae} that there is an $x_0\in
X$ so that the Lipschitz equivalence has a G\^ateaux derivative
$D_f(x_0)$ at $x_0$ and that $f^{-1}$ has a G\^ateaux derivative
at $f(x_0)$. This implies that $D_f(x_0)$ is a surjective isomorphism.

One difficulty in studying questions like  Problem~\ref{pbm3} stems
from the fact that
$\gamma\in\Gamma(X)$ does not necessarily imply that
$f\circ\gamma\in\Gamma(Y)$. However, for spaces with the RNP
this problem is easy to overcome since
$\Gamma$-null sets can be equivalently defined in the
following way. For $\gamma:T\to X$ denote
the Lipschitz constant of $\gamma$ along the $i$-th coordinate by
$$\Lip_i(\gamma)=\sup_{r\in\R,t,t+re_i\in T}
\|\gamma(t+re_i)-\gamma(t)\|/|r| $$
where $e_i$ is the $i$-th unit vector in $\ell_\infty$.
Consider the Fr\'echet space of continuous mappings
$\gamma: T\to X$ for which $\Lip_i(\gamma)<\infty$ for all $i$
equipped with the topology generated by the semi-norms
$\sup_{t\in T}\|\gamma(t)\|$ and $\Lip_k(\gamma)$, and
define
$\tilde\Gamma(X)$ as the closure
of the set of those mappings which depend only
on finitely many coordinates.
Then a Borel set $N\subset X$ is $\Gamma$-null if and only if
$\mu\{t:\gamma(t)\in N\}=0$ for $\gamma$ in a residual set of
$\tilde\Gamma(X)$.

If $f$ is a Lipschitz equivalence from $X$ onto $Y$ then the map
$\varphi_f:\tilde\Gamma(X)\to\tilde\Gamma(Y)$ defined by
$\varphi_f(\gamma) = f\circ\gamma$ is one-to-one and onto.
The map $\varphi_f$ is however highly discontinuous. For giving a
positive answer to Problem~\ref{pbm3} it would be enough to know that
$\varphi_f$ maps residual sets onto residual sets. This however is
false as the following example due to M.~Cs\"ornyei shows:
Let $g:\R\to\R$ be a Lipschitz equivalence such that
$g(0)=-1/4$,
$g'(r)=1$ for $1/4<r<3/4$ and $g'(r)=2$ for $r<1/4$ and $r>3/4$.
Let $\gamma_0\in\tilde\Gamma(\R)$ be defined by
$\gamma_0(t)=t_1$ and let
$$U=\{\gamma\in\tilde\Gamma(\R): \sup_{t\in T} |\gamma(t)-\gamma_0(t)|<1/4,
\Lip_1(\gamma-\gamma_0)<1/4\}.$$
Whenever $\gamma\in U$, we denote $M_\gamma = (g\circ\gamma)^{-1}(1/4,3/4)$,
note that 
$$\int_{M_\gamma} D_1(g\circ\gamma)d\mu= 1/2$$
and that up to a set of $\mu$ measure zero,
$$M_\gamma = \{t\in T: D_1(g\circ\gamma)(t)<5/4\}
=\{t\in T: D_1(g\circ\gamma)(t)<3/2\}.$$
Clearly, for $0<\varepsilon<1/4$
no $\tilde\gamma$ sufficiently close
to $g\circ\gamma(t) +\varepsilon t_1$
can have these properties. Hence
$\varphi_g:\tilde\Gamma(\R)\to\tilde\Gamma(\R)$
maps $U$ onto a nowhere dense subset of $\tilde\Gamma(\R)$.
It follows that $\varphi_{g^{-1}}$ maps the residual set
$\tilde\Gamma(\R)\setminus\varphi_g(U)$ into
$\tilde\Gamma(\R)\setminus U$ which is not residual.

\vglue6pt {\it Acknowledgements}.
We wish to express our thanks to L.~Zaj\'{\i}\v{c}ek for a
number of useful comments and to W.~B.~Johnson and G.~Schechtman
for pointing out, among
other things, that our methods are applicable to the Tsirelson space.

\end{document}